\documentclass[onefignum,onetabnum]{siamart171218}

\usepackage{multirow}
\usepackage{amsfonts,amsmath,amssymb}
\usepackage{graphicx}
 \usepackage{color}
 \usepackage{soulutf8}
\usepackage{units}
\usepackage{graphicx}
\usepackage{url}
\usepackage{textcomp}
\usepackage{array}
\ifpdf
\DeclareGraphicsExtensions{.eps,.pdf,.png,.jpg}
\else
\DeclareGraphicsExtensions{.eps}
\fi

\newsiamremark{remark}{Remark}
\newsiamremark{hypothesis}{Hypothesis}
\crefname{hypothesis}{Hypothesis}{Hypotheses}
\newsiamthm{claim}{Claim}
\newsiamthm{preposition}{Preposition}

\newcommand{\newsiamthmrep}[2]{
	\theoremstyle{plain}
	\theoremheaderfont{\normalfont\sc}
	\theorembodyfont{\normalfont\itshape}
	\theoremseparator{.}
	\theoremsymbol{}
	\newtheorem*{#1}[theorem]{#2}
}

\headers{Node generation on parametric surfaces}{U.\ Duh, G.\ Kosec and J.\ Slak}

\title{Fast variable density node generation on parametric surfaces with application to mesh-free methods
  \thanks{\funding{This work was supported by FWO grant G018916N, the ARRS research core funding no.~P2-0095 and Young Researcher program PR-08346.}}}

\author{Urban Duh\footnotemark[4] \thanks{``Jožef Stefan'' Institute, Department E6, Parallel and
Distributed Systems Laboratory, Jamova cesta 39, 1000 Ljubljana, Slovenia
(\email{urban.duh@student.fmf.uni-lj.si})}
  \and Gregor Kosec\thanks{``Jožef Stefan'' Institute, Department E6, Parallel and Distributed Systems Laboratory, Jamova cesta 39, 1000 Ljubljana, Slovenia
    (\email{gregor.kosec@ijs.si}, \url{http://e6.ijs.si/\string~gkosec/})}
  \and Jure Slak\footnotemark[3] \thanks{Faculty of Mathematics and Physics, University of
  Ljubljana, Jadranska 19, 1000 Ljubljana, Slovenia
      (\email{jure.slak@ijs.si}, \url{http://e6.ijs.si/\string~jslak/}).}
}

\usepackage{amsopn}

\ifpdf
\hypersetup{
  pdftitle={Fast variable density node generation on parametric surfaces with application to mesh-free methods},
  pdfauthor={U. Duh, G. Kosec and J. Slak}
}
\fi

\usepackage{algpseudocode}  % za psevdokodo
\usepackage{algorithm}      % za
\algnewcommand\algorithmicto{\textbf{to}}
\algnewcommand\algorithmicin{\textbf{in}}
\algnewcommand\algorithmicforeach{\textbf{for each}}
\algrenewtext{For}[3]{\algorithmicfor\ #1 $\gets$ #2\ \algorithmicto\ #3\ \algorithmicdo}
\algdef{S}[FOR]{ForEach}[2]{\algorithmicforeach\ #1\ \algorithmicin\ #2\ \algorithmicdo}

\newcommand{\R}{\mathbb{R}}

\newcommand{\T}{\mathsf{T}}
\renewcommand{\b}{\boldsymbol}

\renewcommand{\phi}{\varphi}
\newcommand{\p}{\b{p}}
\newcommand{\xib}{\b{\xi}}
\newcommand{\etab}{\b{\eta}}
\newcommand{\Rem}{\b{R}}

\renewcommand{\r}{\b{r}}

\newcommand{\s}{\vec{s}}
\newcommand{\cb}{\b{c}}
\newcommand{\n}{\b{n}}
\newcommand{\hh}{\hat{h}}

\newcommand{\X}{\mathcal{X}}

\DeclareMathOperator{\obb}{obb}

\begin{document}

\maketitle

% REQUIRED
\begin{abstract}
Domain discretization is considered a dominant part of solution procedures for
solving partial differential equations. It is widely accepted that mesh
generation is among the most cumbersome parts of the FEM analysis and often
requires human assistance, especially in complex 3D geometries.  When using
alternative mesh-free approaches, the problem of mesh generation is simplified
to the problem of positioning nodes, a much simpler task, though still not
trivial. In this paper we present an algorithm for generation of nodes on
arbitrary $d$-dimensional surfaces.  This algorithm complements a recently
published algorithm for generation of nodes in domain interiors, and represents
another step towards a fully automated dimension-independent solution procedure
for solving partial differential equations. The proposed algorithm generates
nodes with variable density on surfaces parameterized over arbitrary parametric
domains in a dimension-independent way in $O(N\log N)$ time. It is also compared
with existing algorithms for generation of surface nodes for mesh-free methods
in terms of quality and execution time.
%It is demonstrated that nodes
%generated with the proposed algorithm can be used in RBF-FD analysis of
%partial differential equations by solving the Poisson equation with mixed
%boundary conditions in irregular 2D and 3D domains.
\end{abstract}

% REQUIRED
\begin{keywords}
Node generation algorithms,
variable density discretizations,
meshless methods for PDEs,
RBF-FD,
Poisson disk sampling
\end{keywords}

% REQUIRED
\begin{AMS}
  65D99, 65N99, 65Y20, 68Q25
\end{AMS}

\section{Introduction}
Generation of nodes on surfaces and their enclosed volumes has many application
in different fields of science and engineering, ranging from computer graphics,
particularly rendering~\cite{bridson2007fast}, to mesh-free numerical analysis
of partial differential equations (PDEs)~\cite{slak2019generation}. In general,
specific applications require specific properties of generated nodes, e.g.\ for
dithering in computer graphics a blue noise distribution is often
desired~\cite{bridson2007fast}, while in mesh-free numerical analysis nodes have
to be positioned regularly enough to support stable numerical approximation of
differential operators~\cite{slak2019generation}.  In the context of mesh-free
analysis positioning algorithms have been developed and tested with different
mesh-free numerical
methods~\cite{fornberg2015fast,kosec2018local,shankar2018robust,slak2019generation,zamolo2018two}.
However, the treatment of boundaries, i.e.\ discretization of surfaces, has
often been overlooked, obtained ad-hoc, or with algorithms for generation of
surface meshes, which are needed in mesh based methods.  Such approach is
conceptually flawed as the whole point of mesh-free methods is to completely
remove  meshing from the procedure. Furthermore, it is also computationally
inefficient, as surface mesh generation algorithms, such
as~\cite{shimada1998automatic} and~\cite{lau1997generation}, spend a great deal
of time generating connectivity relations, only for those relations to be
discarded later. This inefficiencies and increased demand for node generation on
surfaces encourage researches to developed specialized algorithms.

Existing algorithms for point generation on parametric surfaces use different
techniques, which are often generalizations of algorithms for spatial node
generation.  The most straightforward is the naive sampling, which generates
parametric points in the parametric space, and then maps the points to the
surface without any additional processing. Such approach inherently results in a
non-uniform and distorted distribution. To avoid distortions, conformal mappings
can be a promising option~\cite{fornberg2015fast}, if they are easily available
for the surface and if the node placing algorithm in the parametric space
supports variable spacing. Original spatial density can be scaled by square root
of the Gramian of the parametric map to account for the non-uniformity caused by
the map. If not analytically available, conformal mappings can be computed
numerically~\cite{gu2012numerical}, but the computation is expensive and not
worth it, in general.  For specific surfaces, such as spheres or tori, there are
simpler specialized algorithms for point generation~\cite{hardin2016comparison}.
For general surfaces, probabilistic approaches are
available~\cite{diaconis2013sampling,kopytov2015uniform} which generate
uniformly distributed points, but they are often not suitable as node generators
for PDE discretization due to the potentially high irregularity. However, both
naive and probabilistic approaches can be useful to generate initial
distributions for iterative discretization generation schemes, such as minimal
energy nodes~\cite{hardin2004discretizing}, energy functional minimization~\cite{vlasiuk2018fast}
 or via dynamic simulation with
attractive/repulsive forces~\cite{shimada1998automatic, kosec2018local}, which
are also commonly used in graphics community for various purposes, such as
texture generation~\cite{turk2001texture}.  Another approach to surface node
generation was presented in the paper by Shankar, Kirby and
Fogelson~\cite{shankar2018robust}, which presents surface reconstruction,
surface node generation and spatial node generation algorithms.
The algorithms were used to obtain discretizations suitable for strong form
mesh-free methods, and the surface node generation technique they used is
called \emph{supersampling-decimation}, which samples the parametric space with
increased density, maps the points and then decimates the mapped points to
conform to the required nodal spacing.

Existing algorithms for node generation on curved surfaces have their
shortcomings. Algorithms based on conformal maps are practical only for specific
classes of surfaces, probabilistic algorithms usually do not produce
distributions of sufficient quality and iterative schemes are needed to improve
them, making them less efficient. The naive and supersampling-decimation
approaches have their benefits in simplicity and speed in some cases, but the
published version in~\cite{shankar2018robust} only deals with cases where the
surface is homeomorphic to a sphere $\mathbb{S}^1$ or $\mathbb{S}^2$, parametric
domain is a rectangle and nodal spacing is constant.  In this paper we present
an algorithm for discretization of surfaces that works in arbitrary dimensions
with variable nodal spacing, with irregular surfaces and parametric domains, but
at the cost of requesting the user to supply the Jacobian $\nabla \r$ of the
surface map.  It has guaranteed $O(N \log N)$ time complexity regardless of the
domain shape and nodal spacing and the running time is comparable to published
algorithms

The rest of the paper is organized as follows. The proposed surface node placing
algorithm is presented in~\cref{sec:alg} along with possible generalizations,
the comparison with existing surface placing algorithms is presented
in~\cref{sec:comp}, its use in meshless numerical simulations is presented
in~\cref{sec:num} and conclusions are presented in~\cref{sec:conc}.

\section{Node generation algorithm}
\label{sec:alg}
Boundaries of computational domains can be represented in different forms.  Most
common ones are as parametrizations (possibly split into patches), such as
produced with non-uniform rational B-splines (NURBS)~\cite{piegl2012nurbs} or
Radial Basis Functions (RBFs)~\cite{carr2003smooth}, as
level-sets~\cite{zhao2001fast} or as subdivision
surfaces~\cite{litke2001fitting}.  Different representations have different
desirable and undesirable properties, as discussed in
e.g.~\cite{stroud2006boundary}.  We will assume that a parametric
representation of the boundary in question is given as
$\r \colon \Xi \subset \R^{d_\Xi} \to \partial \Omega \subset \R^d$,
along with the Jacobian $\nabla \r$. Most of the algorithms described in the
introduction also work with parametric representations, and a way to construct
them is also offered in~\cite{shankar2018robust}. Some surface quantities
obtained from higher order derivatives, such as curvature, may also be required
to solve PDEs, but they will not be used during node generation.
The elements of \emph{parametric space} $\Xi$ will be called \emph{parameters},
and the elements of \emph{target space} $\partial \Omega$ will be called
\emph{points} or \emph{nodes}, when specific emphasis on discretization is
desired. Usually, the dimension $d_\Xi$ will be $d-1$, but as described
in~\cref{sec:gen-rem} this is not a requirement.

Given a nodal spacing function $h\colon \partial\Omega \subset \R^d \to (0,
\infty)$, we wish to place nodes on $\partial \Omega$, such that spacing around
node $\p\in\partial\Omega$ is approximately equal to $h(\p)$.

The proposed node placing algorithm takes a regular parametrization $\r$, its
Jacobian $\nabla \r$, a nodal spacing function $h$ and a set of ``seed
parameters'' $\X$ from $\Xi$ as input. If no seed parameters are supplied by the
user, the algorithm chooses a random starting parameter inside $\Xi$.  It
returns a set of regularly distributed nodes on the surface $\partial \Omega$,
conforming to the spacing function $h$.

The node spacing does not take place directly in the target space $\R^d$.
Instead, we place parameters in the parametric space $\Xi$ using the same
principle as for spatial node placing. However, the distance between two
parameters $\xib_1$ and $\xib_2$ is not chosen to be local spacing $h$ but is
instead computed in such a way, that the distance between points $\r(\xib_1)$
and $\r(\xib_2)$ is approximately $h$.  A spatial search structure of points in
$\R^d$ is maintained to check for proximity violations.

% usual description
The node placing algorithm processes nodes sequentially. Initially, seed parameters are put
in a queue and their corresponding points in the spatial search structure.
In each iteration, a parameter $\xib_i$ is dequeued, and
expanded into a set of candidate parameters $H_i$, by generating a set of directions,
represented by unit vectors $\s_{i, j}$ that approximately uniformly cover all possible spatial
directions in $\Xi$.

To derive how far from $\xib_i$ a candidate parameter must be placed to achieve appropriate
distance between points,
we consider a candidate parameter $\etab_{i, j} \in H_i$ in the direction $\s_{i, j}$ from $\xib_i$,
\begin{equation} \etab_{i, j} = \xib_i + \alpha_{i, j} \s_{i, j}, \end{equation}
for some distance $\alpha_{i, j} > 0$. Denote the point corresponding to $\xib_i$ as
$\p_i = \r(\xib_i)$. We would like the distance between the candidate point $\r(\etab_{i, j})$ and
$\p_i$ in $\partial \Omega$ to be equal to $h(\p_i)$
\begin{equation} \label{eq:h-wish}
\|\r(\etab_{i, j}) - \r(\xib_i)\| = h(\p_i),
\end{equation}
where $\|\cdot\|$ is the standard Euclidean distance.
Using first order Taylor's expansion, we write
\begin{equation} \label{eq:tay-lin}
\r(\etab_{i, j}) = \r(\xib_i + \alpha_{i, j} \s_{i, j}) \approx \r(\xib_i) + \alpha_{i, j} \nabla \r (\xib_i) \s_{i, j}.
\end{equation}
Substituting the linear approximation~\eqref{eq:tay-lin}
in~\eqref{eq:h-wish}, we obtain an equation for $\alpha_{i, j}$ as
\begin{equation} \label{eq:alpha}
  h(\p_i) = \|\r(\xib_i) + \alpha_{i, j} \nabla \r (\xib_i) \s - \r(\xib_i)\| = \alpha_{i, j} \| \nabla \r (\xib_i) \s_{i, j} \|,
\end{equation}
where we took into account the positivity of $\alpha_{i, j}$. This
equation is solved for $\alpha_{i, j}$ to obtain the candidate parameter
\begin{equation}
  \etab_{i, j} = \xib_i + \frac{h(\p_i)}{ \| \nabla \r (\xib_i) \s_{i, j} \|} \s_{i, j}, \quad \alpha_{i, j} =  \frac{h(\p_i)}{ \|
  \nabla \r (\xib_i) \s_{i, j} \|}.
\end{equation}
Higher order terms of the approximation could be used in~\eqref{eq:tay-lin}, but the
equation~\eqref{eq:alpha} would also be of higher order and higher derivatives of $\r$
would be required.

The set of candidates $H_i$ is generated from directions $\vec{s}_{i,j}$ as
\begin{equation}
  H_i = \left\{ \xib_i + \frac{h(\p_i)}{ \| \nabla \r (\xib_i) \s_{i, j} \|} \s_{i, j} \mid \s_{i,
  j} \in \text{directions}(\xib_i) \right\}.
\end{equation}
Candidate parameters that lie outside of the parametric space $\Xi$ are discarded.
Additionally, candidate parameters whose corresponding points lie too close to already accepted
points are also rejected. The remaining candidate parameters are enqueued for expansion and their
corresponding points are added to the spatial search structure.

When checking the distance from a candidate point to its closest already existing point,
the algorithm compares the found distance with $\|\r(\etab_{i, j}) - \r(\xib_i)\| =: \hh_{i, j}$
instead of $h(\p_i)$. Unless this adjustment is made, no candidates would be accepted in areas
where $\hh_{i, j}$ is smaller than $h(\p)$, since point $\p_i$ itself would be too close.
On curves, this might even cause the algorithm to terminate prematurely.

\begin{figure}[h]
  \centering
  \includegraphics[width=0.4\linewidth]{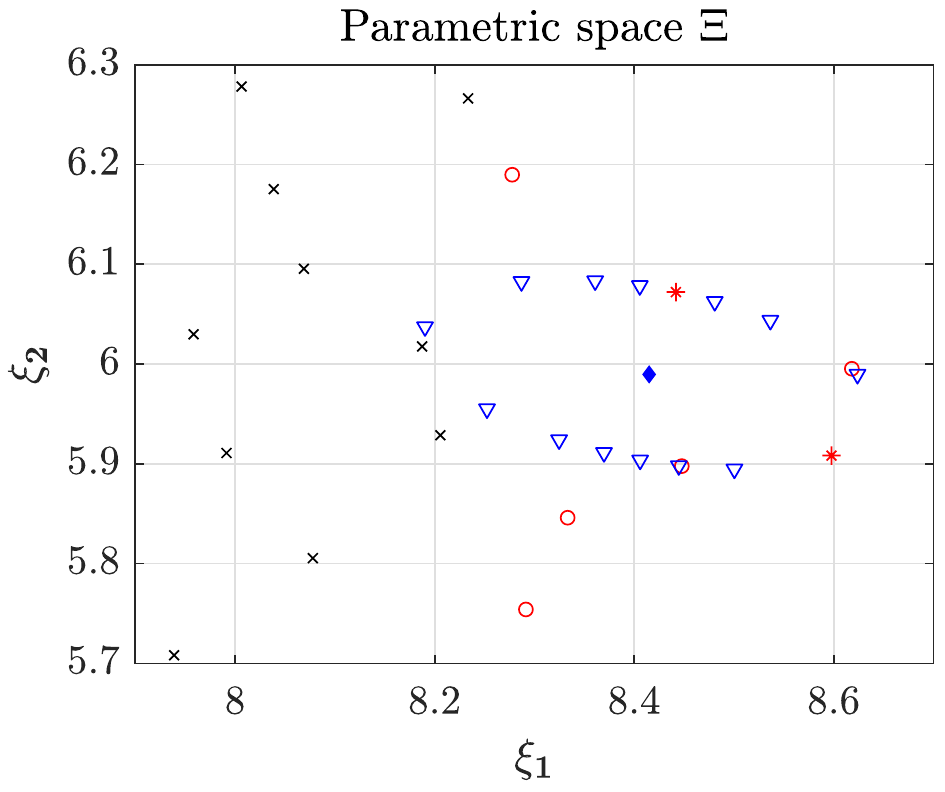}
  \includegraphics[width=0.59\linewidth]{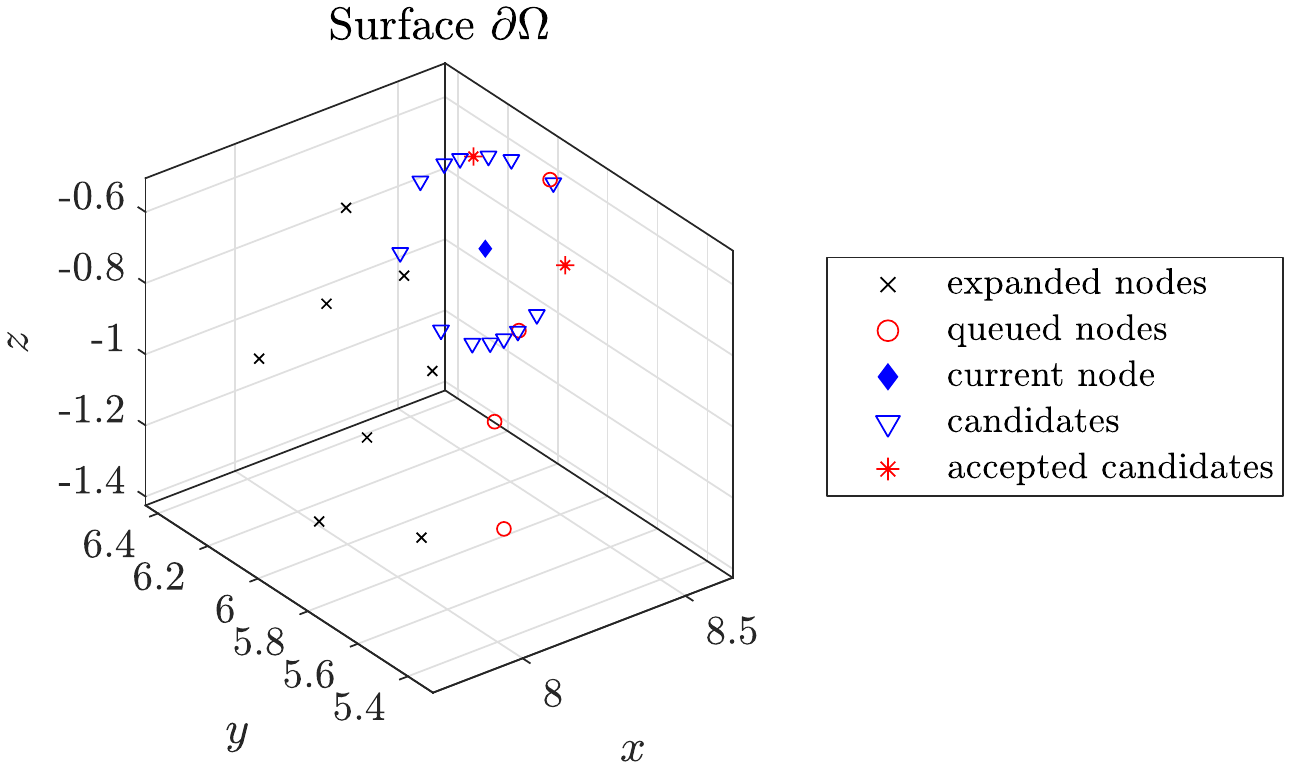}

  \caption{Illustration of candidate generation by the proposed algorithm in
  parametric space $\Xi$ (left) and main
  domain $\partial \Omega$ (top), around parameter $\xib = (8.42, 5.99)$ with spacing $h = 0.23$. }
  \label{fig:gen-ill}
\end{figure}

\Cref{fig:gen-ill} illustrates the process of candidate generation on a surface parameterized with
$\r(\xi_1, \xi_2) = (\xi_1, \xi_2, 3 \sin \xi_1 \sin \xi_2)$.
Parameters in parametric space $\Xi$ (left) are generated  in a way that when mapped
to the main domain $\partial \Omega$ (right), they are approximately $h$ apart.

\begin{figure}[h!]
  \centering
  \includegraphics[width=0.19\linewidth]{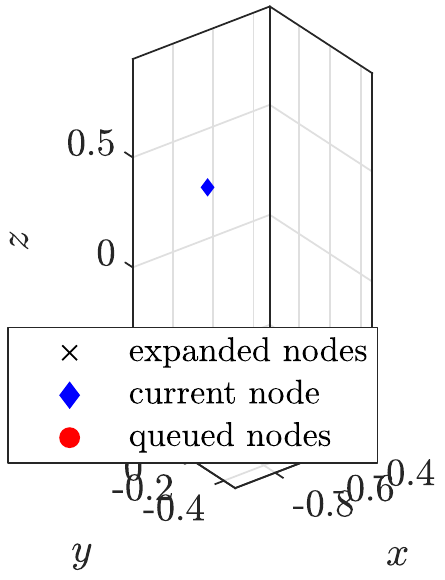}
  \includegraphics[width=0.19\linewidth]{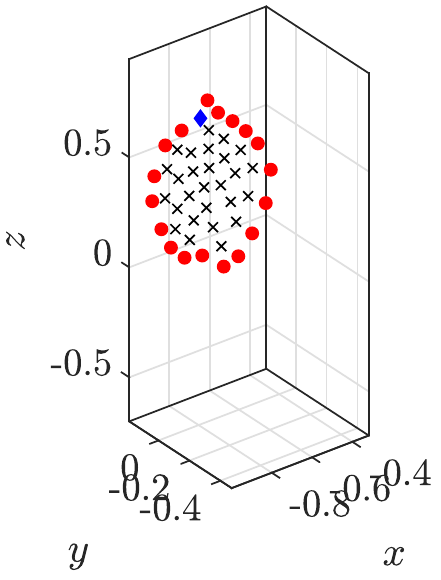}
  \includegraphics[width=0.19\linewidth]{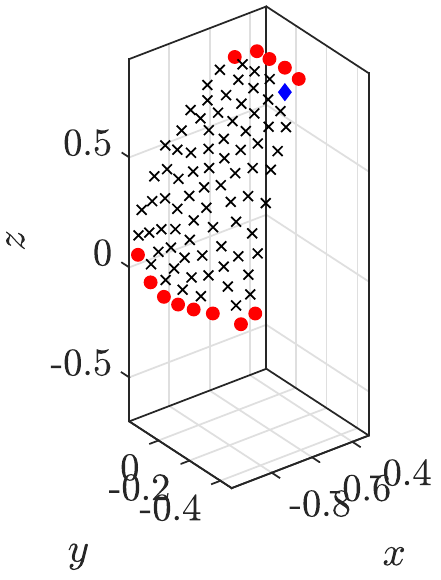}
  \includegraphics[width=0.19\linewidth]{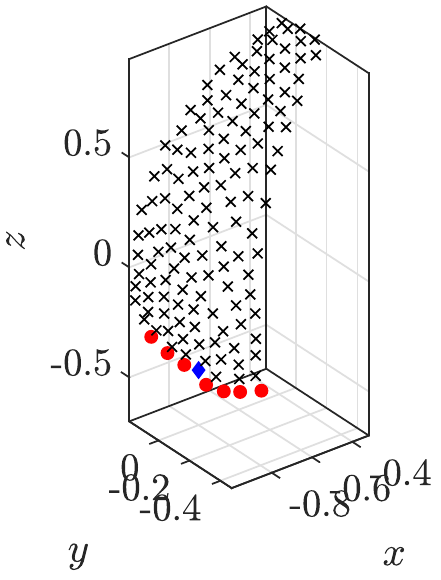}
  \includegraphics[width=0.19\linewidth]{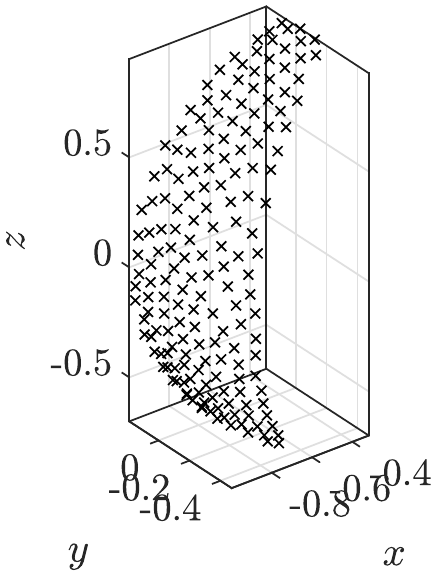}

  \includegraphics[width=0.2\linewidth]{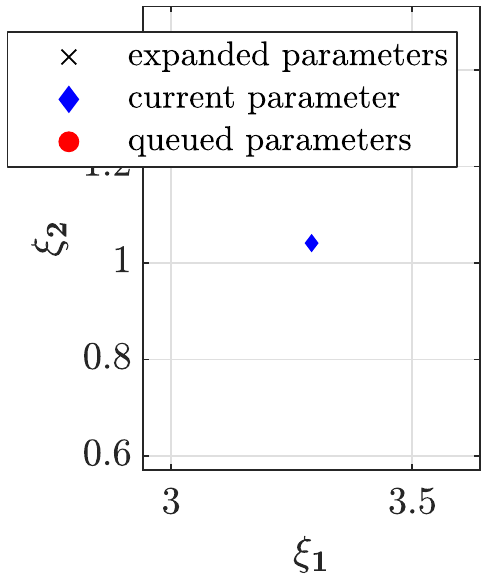}
  \includegraphics[width=0.19\linewidth]{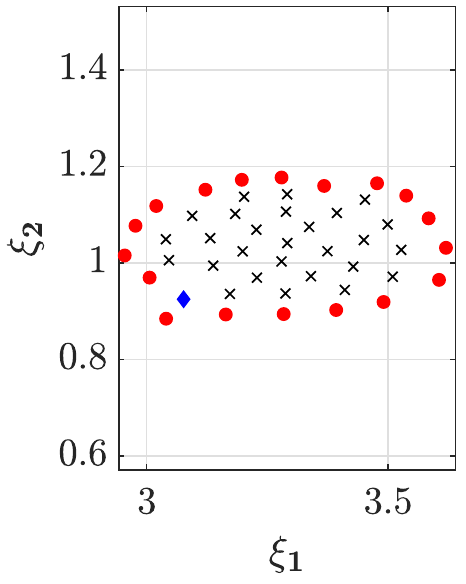}
  \includegraphics[width=0.19\linewidth]{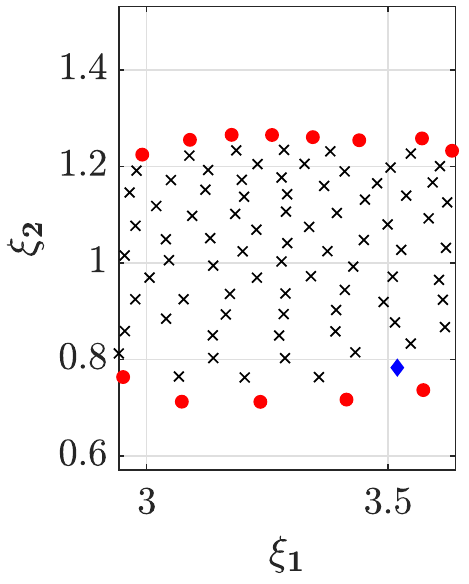}
  \includegraphics[width=0.19\linewidth]{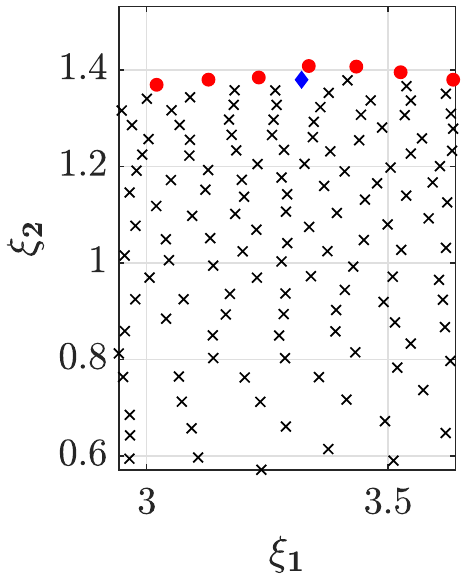}
  \includegraphics[width=0.19\linewidth]{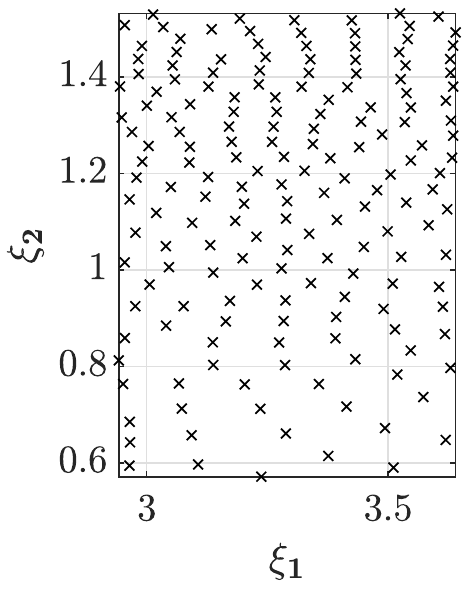}

  \caption{Execution illustration of the proposed algorithm algorithm (left to
  right) in parametric domain $\Xi$ (bottom) and main domain $\partial \Omega$ (top).
  Part of a unit sphere was sampled with nodal spacing $h = 0.08$.}
  \label{fig:exe-ill-run}
\end{figure}

\Cref{fig:exe-ill-run} illustrates the execution of the proposed algorithm on a
part of a unit sphere (non-standardly) parameterized by
$\r(\xi_1, \xi_2) = (\cos \xi_1 \sin \xi_2^2, \sin \xi_1 \sin \xi_2^2, \cos \xi_2^2)$ with a constant nodal spacing.
The algorithm begins with a single node (leftmost image) and then expands it in all directions.
Subsequent images show progress of the proposed algorithm and the final result
(rightmost image).

An efficient implementation with an implicit queue contained in the array of final points and
the $k$-d tree spatial structure~\cite{moore1991intoductory} is presented as~\cref{alg:pa}.

\begin{algorithm}[h!]
  \small
  \caption{Proposed surface node placing algorithm.}
  \label{alg:pa}
  \vspace{5pt}
%  \textbf{Input:} Domain $\Omega$ and its dimension $d$. \\
  \textbf{Input:} Parametric domain $\Xi \subseteq \R^{d_\Xi}$, given by its characteristic
  function $\chi_\Xi\colon\R^{d_\Xi} \to \{0,1\}$. \\
  \textbf{Input:} Parametrization $\r \colon \Xi \to \partial \Omega \subset \R^d$ and its Jacobian
  matrix $\nabla \r$. \\
  \textbf{Input:} A list of starting parameters $\X$ from the parametric domain $\Xi$. \\
  \textbf{Input:} A nodal spacing function $h\colon \partial \Omega \subset \R^d \to (0, \infty)$. \\
  \textbf{Input:} Number of candidates generated in each iteration $n$. \\
  \textbf{Output:} A list of points in $\partial \Omega$ distributed according to spacing function $h$.
  \begin{algorithmic}[1]
    \Function{proposedAlgorithm}{$\Xi, \r, \nabla \r, h, \X, n$}
    \State $T \gets \Call{kdtreeInit}{\r(\X)}$  \Comment{Initialize spatial search structure on points $\r(\X)$.}
    \State $i \gets 0$ \Comment{Current node index.}
    \While{$i < |\X|$} \Comment{Until the queue is not empty.} \label{ln:while}
    \State $\xib_i \gets \X[i]$ \Comment{Get next parameter values from the start of the queue.}
    \State $\p_i \gets \r(\xib_i)$ \Comment{Compute the corresponding node in $\partial \Omega$.}
    \State $h_i \gets h(\p_i)$ \Comment{Compute its nodal spacing.}
    \ForEach{$\s_{i, j}$}{\Call{candidates}{$n$}} \Comment{Loop through random unit vectors.} \label{ln:candidates}
    \State $\etab_{i, j} \gets \xib_i + \frac{h_i}{ \| \nabla \r (\xib_i) \s_{i, j} \|} \s_{i, j}$ \Comment{Calculate new candidate.}
    \If{$\etab_{i, j} \in \Xi$} \Comment{Discard candidates outside the parametric domain.}
    \State $\cb_{i, j} \gets \r(\etab_{i, j})$ \Comment{Compute the candidate point in $\partial \Omega$.}
    \State $\hh_{i, j} \gets \|\cb_{i, j} - \p_{i, j}\|$ \Comment{Compute the actual spacing.}
    \State $\n_{i, j} \gets \Call{kdtreeClosest}{T, \cb_{i, j}}$ \Comment{Find nearest node for proximity test.}
    \If{$\|\cb_{i, j} - \n_{i, j}\| \geq \hh_{i, j}$} \Comment{Test that $\cb_{i, j}$ is not too close to other nodes.}
    \State \Call{append}{$\X, \etab_{i,j}$}  \Comment{Enqueue $\etab_{i,j}$ as the last element of $\X$.}
    \State \Call{kdtreeInsert}{$T, \cb_{i,j}$}  \Comment{Insert $\cb_{i,j}$ into the spatial search structure.}
    \EndIf
    \EndIf
    \EndFor
    \State $i \gets i + 1$ \Comment{Dequeue current parameter and move to the next one.}
    \EndWhile
    \State \Return $\r(\X)$
    \EndFunction
  \end{algorithmic}
\end{algorithm}

The proposed algorithm includes generation of random unit direction vectors on
line~\ref{ln:candidates},
that needs to be clarified.
There are different ways of generating unit vectors that cover all possible directions well,
according to discussion in~\cite{slak2019generation}, \textit{randomized pattern candidates}
technique shall be used. A set of random unit vectors is obtained by
randomly rotating a fixed discretization of a $d$-dimensional unit ball.
The set of candidates $\textsc{candidates}(n)$ on a unit ball in 2D is obtained simply by
\begin{equation}
\textsc{candidates}(n) = \left\{ (\cos\phi, \sin\phi); \ \phi \in \{0, \phi_0, 2\phi_0, \ldots,
(n-1)\phi_0\}, \phi_0 = \frac{2\pi}{n} \right\}.
\end{equation}
In $d$ dimensions, the discretization of a unit ball is obtained by recursively
discretizing appropriate $d-1$ dimensional slices along the last coordinate.

The number of generated candidates $n$ in each iteration on line~\ref{ln:candidates} is also a free
parameter in our implementation of the proposed algorithm and recommendations based on the dimension
$d_\Xi$ are given in~\cite{slak2019generation}. We will use $n = 2$ when $d_\Xi = 1$ and $n = 15$ when $d_\Xi = 2$ in our analyses,
unless stated otherwise.

\subsection{Time complexity analysis}
\label{sec:time-theoretical}
We will derive the time complexity in terms of the number of generated nodes.
Let us denote the number of starting points with $N_s$, the number of final points with $N$,
the cost of spatial search structure precomputation, query and insertion with
$P(N_s), Q(N)$ and $I(N)$ respectively, the cost of evaluating $\r$ with $e_1$ and the cost of
evaluating $\nabla \r$ with $e_2$. All other operations are assumed to have (amortized) constant
cost.

The main while loop (line~\ref{ln:while}) iterates exactly $N$ times and the inner for loop
(line~\ref{ln:candidates}) iterates $n$ times. Each iteration of the inner for loop executes one query and one $\r$
evaluation. Every candidate was inserted once and $\nabla \r$ is evaluated once for each node when expanding it, since
its value can be stored for later use in the inner for loop.
Therefore, the total time complexity of the proposed algorithm is equal to
\begin{align}
  T_\text{PA-general} &= P(N_s) + \sum_{i = 1}^{N}(I(N) + e_1 + e_2 + \sum_{j = 1}^{n}(e_1+ Q(N))) \label{eq:time-long} \\
                &= O(P(N_s) + N(n+1)e_1 + N e_2 + NnQ(N) + NI(N)).  \nonumber
\end{align}

Since we are using a $k$-d tree as the spatial search structure, we know that $P(N_s) = O(N_s \log N_s), Q(N) = O(\log N)$ and $I(N) = O(\log N)$. Both $e_1$ and $e_2$ are
also usually constant cost operations and $N_s = O(N)$ (usually even $N_s \ll N$). This simplifies
equation~\eqref{eq:time-long} to
\begin{equation}
  T_\text{PA} = O(N_s \log N_s + nN + nN \log N) = O(nN \log N). \label{eq:time-short}
\end{equation}

Other data structures can be used in special cases. If $h$ is constant, a background grid
with spacing $O(h/\sqrt{d})$ can be faster, but less space efficient, supporting
query and insert operations in $O(1)$. Additionally, when
when sampling simple curves from a single starting point, the search structure is unnecessary
as the parameters are sampled in the interval and we can simply advance in both directions
and the only conflict can appear when the two ends of the curve potentially meet.

In this article, we will be using the general $k$-d tree version of the algorithm
unless stated otherwise explicitly.

\subsection{Remarks on generalizations}
\label{sec:gen-rem}

In general, the parametrization function $\r$ does not need to map to a smooth boundary $\partial
\Omega$ of a domain $\Omega$. However, having orientability and co-dimension one allows the
algorithm to also uniquely generate unit normals from $\nabla \r$.
Normal generation aside, neither orientability nor co-dimension are requirements
for the algorithm. The algorithms works in exactly the same manner if
$d_\Xi < d-1$, e.g.\ for a curve embedded in 3D space.

\begin{figure}[h]
  \centering
  \includegraphics[width=0.46\linewidth]{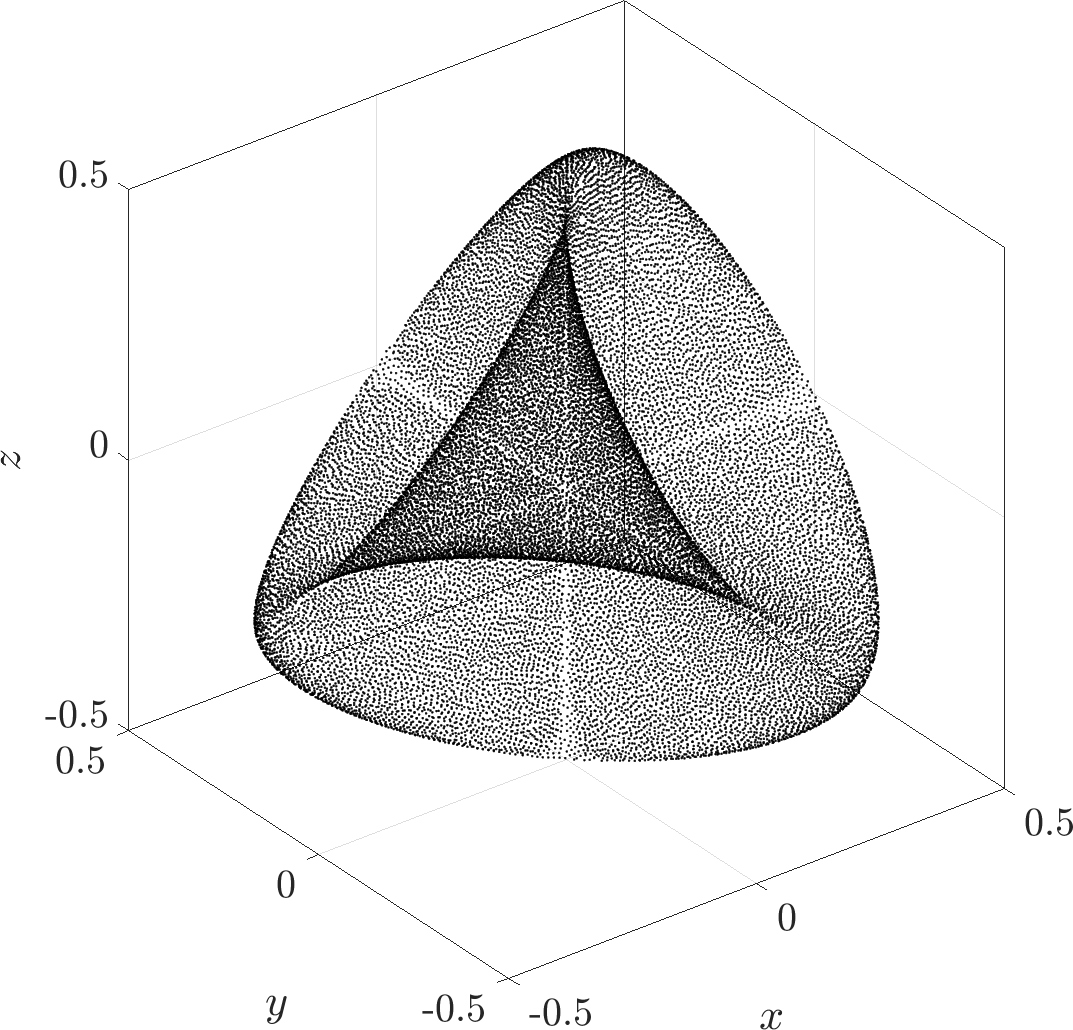}
  \includegraphics[width=0.53\linewidth]{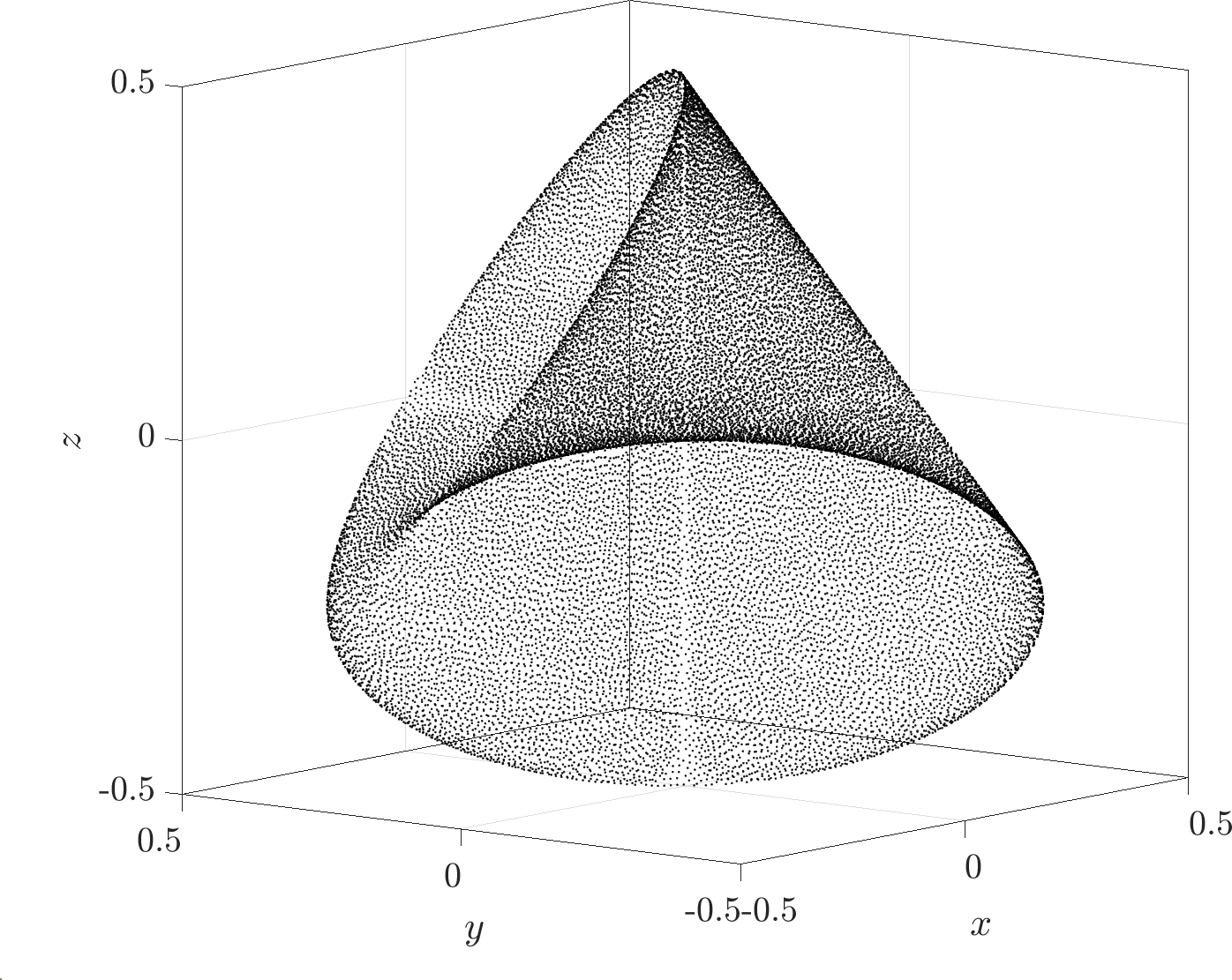}
  \caption{Nodes generated by the proposed algorithm on the Roman
  surface~\eqref{eq:roman} with $R = 1, h = 0.005$.}
  \label{fig:manifold}
\end{figure}

Additionally, there are also no problems with self-intersecting surfaces.
As an example of that, the discretization of the (non-orientable)
Roman surface is shown in~\cref{fig:manifold}.
The Roman surface is a mapping of the real projective plane $\R\textup{P}^2$ in $\R^3$, given by
\begin{equation} \label{eq:roman}
  \r(\theta, \phi) = (
  R^2 \cos \theta \sin \theta \sin \phi,
  R^2 \cos \theta \sin \theta \cos \phi,
  R^2 \cos^2 \theta \cos \phi \sin \phi).
\end{equation}

%Nodes could also be generated on a non self-intersection 4-D
%embedding.

Although all examples in this paper only deal with closed surfaces, the
generalization to bounded surfaces is straightforward. First, the boundary
is to be discretized, using the
proposed algorithm, followed by discretization of remaining of the surface
using the generated boundary nodes as seed nodes.

There is also a possibility to extend the algorithm to surfaces defined by
multiple possibly
intersecting patches, such as models created by Computer aided design (CAD)
software or
parameterized submanifolds. The algorithm can discretize the surface patch by patch and the spatial
search structure keeps all the information about stored nodes from the previous
patches, to check for spacing violations. Problems might arise on patch joints, similarly to
front-joints, which could be dealt with post-process algorithms~\cite{kosec2018local}.
However, further discussion on this topic is out of scope of this paper and is left for future work.

\section{Discussion and comparison}
\label{sec:comp}
The proposed algorithm is compared with the supersampling-decimation technique, recently
published by Shankar, Kirby and Fogelson~\cite{shankar2018robust},
and the naive sampling algorithm,
which was also used in~\cite{shankar2018robust} to show the significance of the supersampling-decimation technique.
A brief description of both algorithms follows, using
the same notation as in~\cref{sec:alg}. When needed, we will refer to the naive algorithm as NA,
to the supersampling-decimation approach as SD and to the proposed algorithm as PA.

\subsection{Existing algorithms}
\subsubsection{Naive parametric sampling}
Naive parametric sampling attempts to generate a discretization of a surface parameterized with
$\r \colon \Xi \subset \R^{d_\Xi} \to \partial \Omega \subset \R^d$,
by discretizing $\Xi$ with nodal spacing $h$, obtaining a set of parameters $\X_\Xi$.
The discretization points are obtained by simply mapping elements of $\X_\Xi$ to the surface,
i.e.\ the resulting set of points $\X$ is given by $\X = \r(\X_\Xi)$.

This type of sampling is useful for its simplicity, especially if $\Xi$ is box-shaped, $h$
constant and gradients $\nabla \r$ are not too large. The algorithm is included in this paper
mostly as a reference to put results in perspective.

\subsubsection{Supersampling-decimation}
The supersampling algorithm is based on the naive algorithm, except that
the discretization of the parametric space $\Xi$ uses spacing $h_\Xi := h/\gamma$ instead of $h$,
where $\gamma > 0$ is called the \emph{supersampling factor}.
This is done in order to generate enough nodes even where the mapping $\r$ might cause them to
spread out on the surface. After the set of parameters $\X_\Xi$ with spacing $h/\gamma$
has been generated, \emph{decimation} or thinning is performed, to keep only the appropriate nodes.
The parameters are mapped to the surface in sequence, accepting only the points
that are not too close to already accepted ones. This requires the use of a spatial search
structure that supports ball queries. The end result is a set of
nodes $\X$ in $\partial \Omega$, that are spaced approximately by $h$.

The algorithm as described in~\cite{shankar2018robust} assumes that $\Xi$ is either a line
or a rectangle and that spacing $h$ is constant. Additionally, $\gamma$ is not chosen directly,
but indirectly by estimating the number of nodes $N$ on the generated surface as
$N = |\partial \Omega| / h^d$ and generating $\tau N$ of them, with spacing
$h_\X = \sqrt[d_\Xi]{|\Xi|/(\tau N)} = h\sqrt[d_\Xi]{|\Xi|/(\tau |\partial \Omega|)}$,
effectively choosing $\gamma = \sqrt[d_\Xi]{\tau |\partial \Omega|/|\Xi|}$.
As $|\partial \Omega|$ is not known directly, is is estimated with the surface of an (oriented)
bounding box of $\Omega$.

We will similarly take into the account the scaling due to different surface areas, but will scale
$h$ directly by $\tau$, using $\gamma = \tau \sqrt[d_\Xi]{|\partial \Omega|/|\Xi|}$.
Value $\tau = 5$ will be sufficient in most cases and will be used unless otherwise specified.

\subsection{Setup for comparison}
In this paper we focus on generating nodes for use in meshless numerical
analysis and therefore all analyses are done with guidelines established
in~\cite{slak2019generation} in mind. These include local regularity, minimal
spacing requirements, computational efficiency and the number of tuning parameters.
In the following discussion we analyze algorithm presented
in~\cref{sec:alg} and compare it to the naive and supersampling algorithms.

In most of the analyses we will use a polar curve used
in~\cite{shankar2018robust},
given by \begin{equation}
r_p(\phi) = |\cos(1.5 \phi)|^{\sin(3 \phi)} \label{eq:polar-curve}, \quad
  \r_p(\phi) = (r_p(\phi) \cos\phi, r_p(\phi) \sin\phi), \quad \phi \in [0, 2\pi)
\end{equation}
and a heart-like surface in 3D, given by
\begin{equation}
\r_h(u, v) = (\sqrt{1 - v^2} \cos(u) + v^2, \sqrt{1 - v^2}\sin(u), v), \quad
(u, v) \in [0, 2 \pi) \times [-1, 1). \label{eq:heart-surf}
\end{equation}
Both of these parametrizations have large variations in absolute value of partial
derivatives, which makes them good candidates for analyses.

All three algorithms were implemented in C++ using the same
library for $k$-d tree spatial search structure (nanoflann~\cite{blanco2014nanoflann})
and the same linear algebra library (Eigen~\cite{eigenweb}) to ensure as fair
comparison as possible. This implementation of the proposed algorithm is included in the Medusa library~\cite{medusa}, a C++ library focused on tools for solving Partial Differential Equations with strong-form meshless methods. Its standalone C++ implementation is also available at~\cite{standalonepa}.

\subsection{Local regularity}
\begin{figure}[h]
  \centering
  \includegraphics[width=0.335\linewidth]{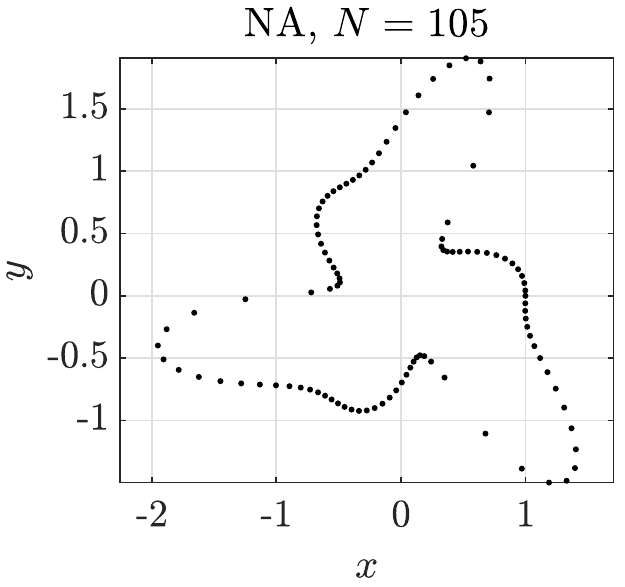}
  \includegraphics[width=0.32\linewidth]{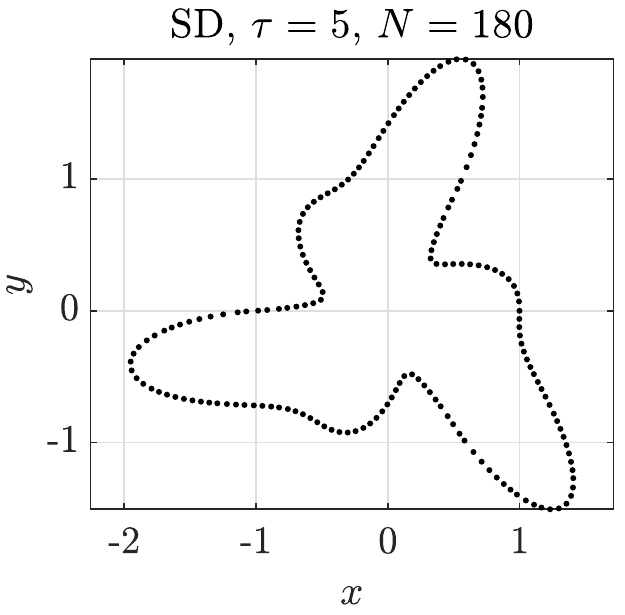}
  \includegraphics[width=0.32\linewidth]{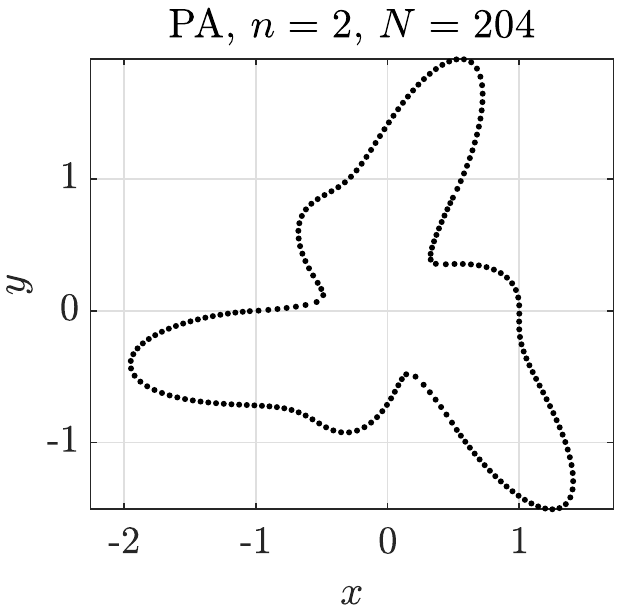}
  \caption{Comparison of different algorithms on a 2D polar curve from~\eqref{eq:polar-curve} sampled with $h = 0.06$}
  \label{fig:nodes-2d}
\end{figure}

We begin our analysis by visually comparing the node sets, generated with different algorithms.
In~\cref{fig:nodes-2d} we can see their performance on the polar curve $\r_p$.
It is clear that the naive algorithm does not perform well on curves with variable derivatives.
We can also see some irregularly big gaps between nodes generated by the supersampling algorithm,
since the supersampling algorithm is based on the naive algorithm. However, this can be improved by
choosing a bigger value of the supersampling parameter $\tau$, at greater cost to execution time
and memory. Nodes generated by the proposed algorithm only have one visually bigger gap, where two
sides of the parametric domain meet.

\begin{figure}[h]
  \centering
  \includegraphics[width=0.32\linewidth]{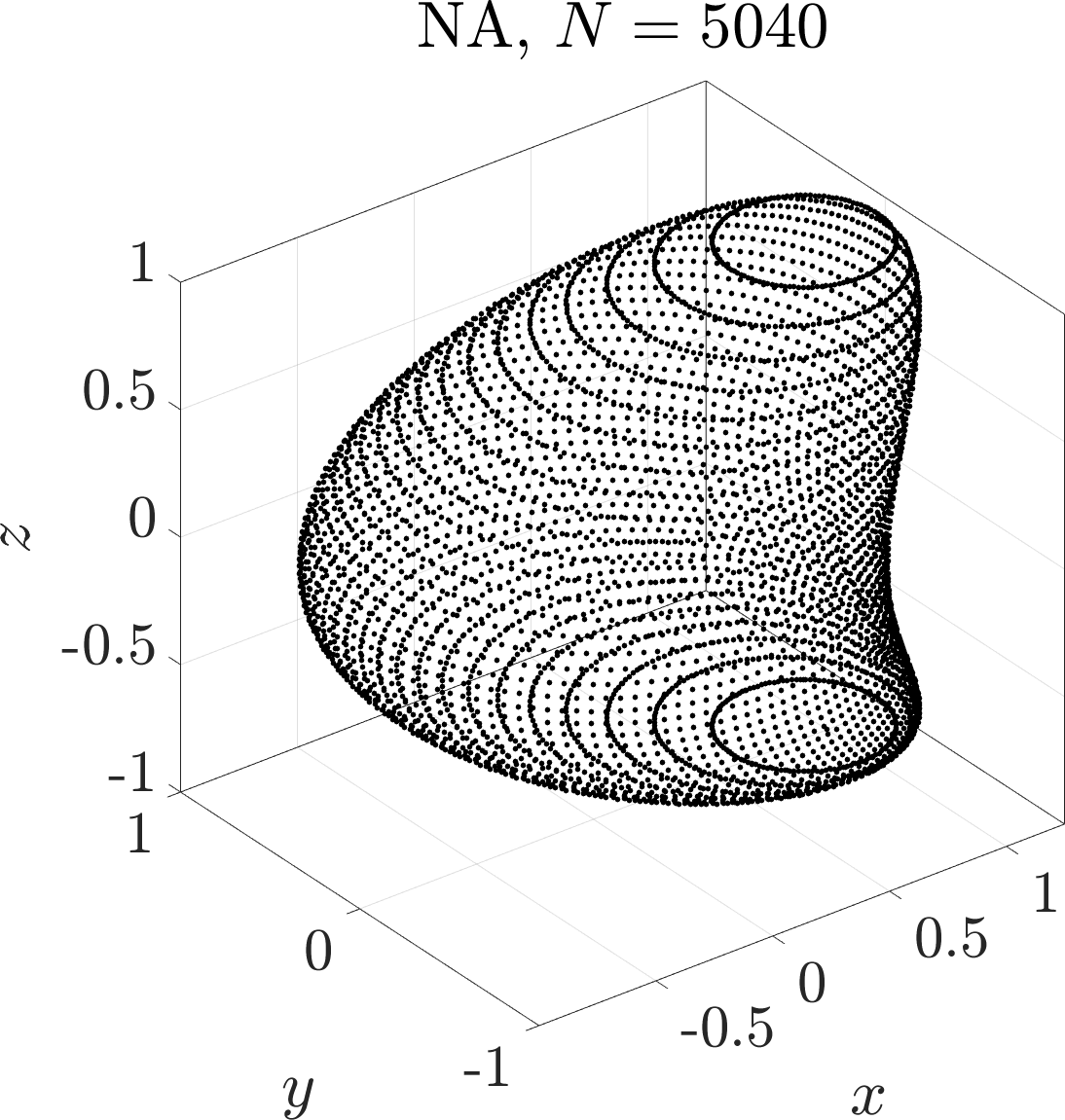}
  \includegraphics[width=0.326\linewidth]{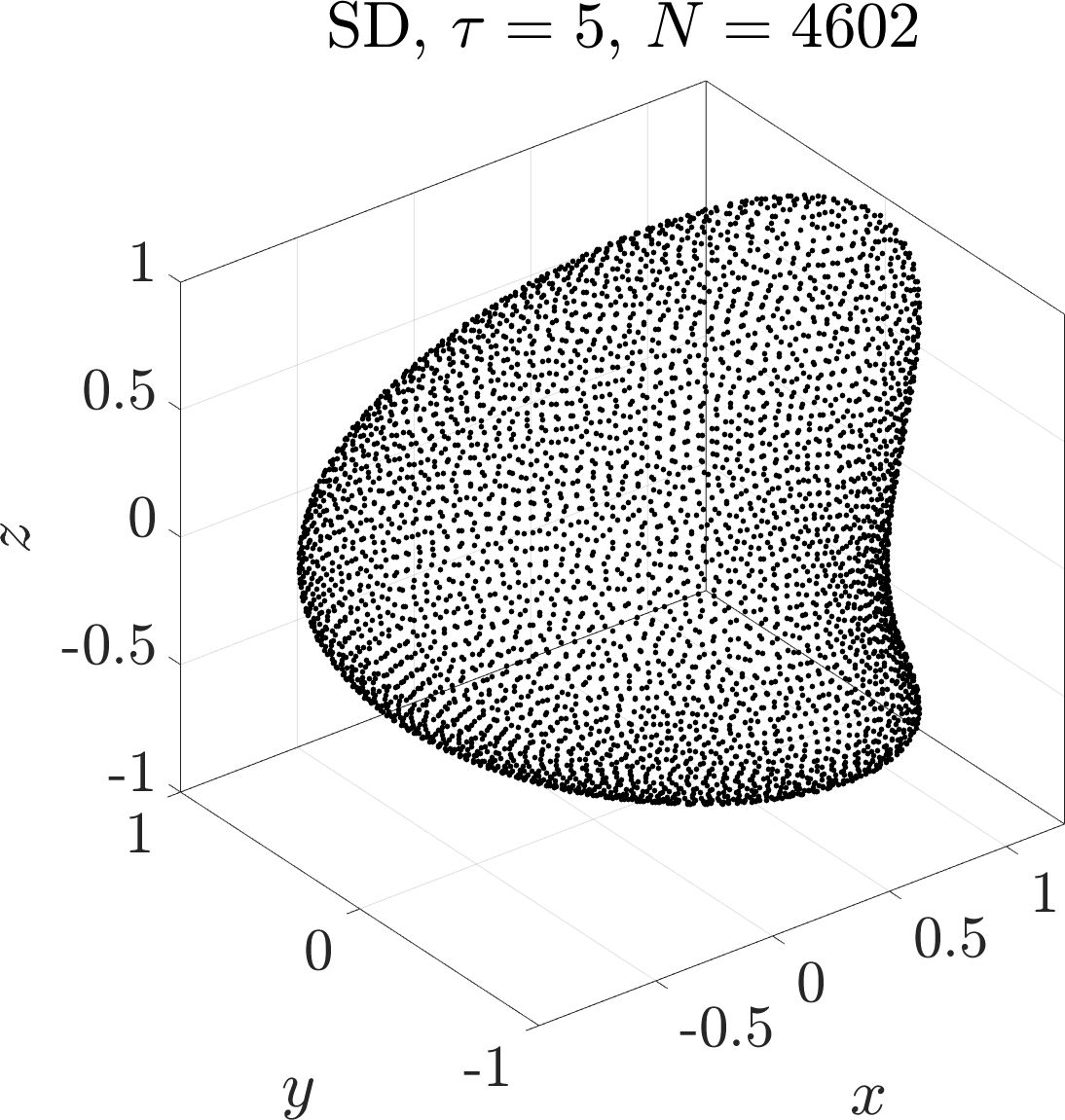}
  \includegraphics[width=0.32\linewidth]{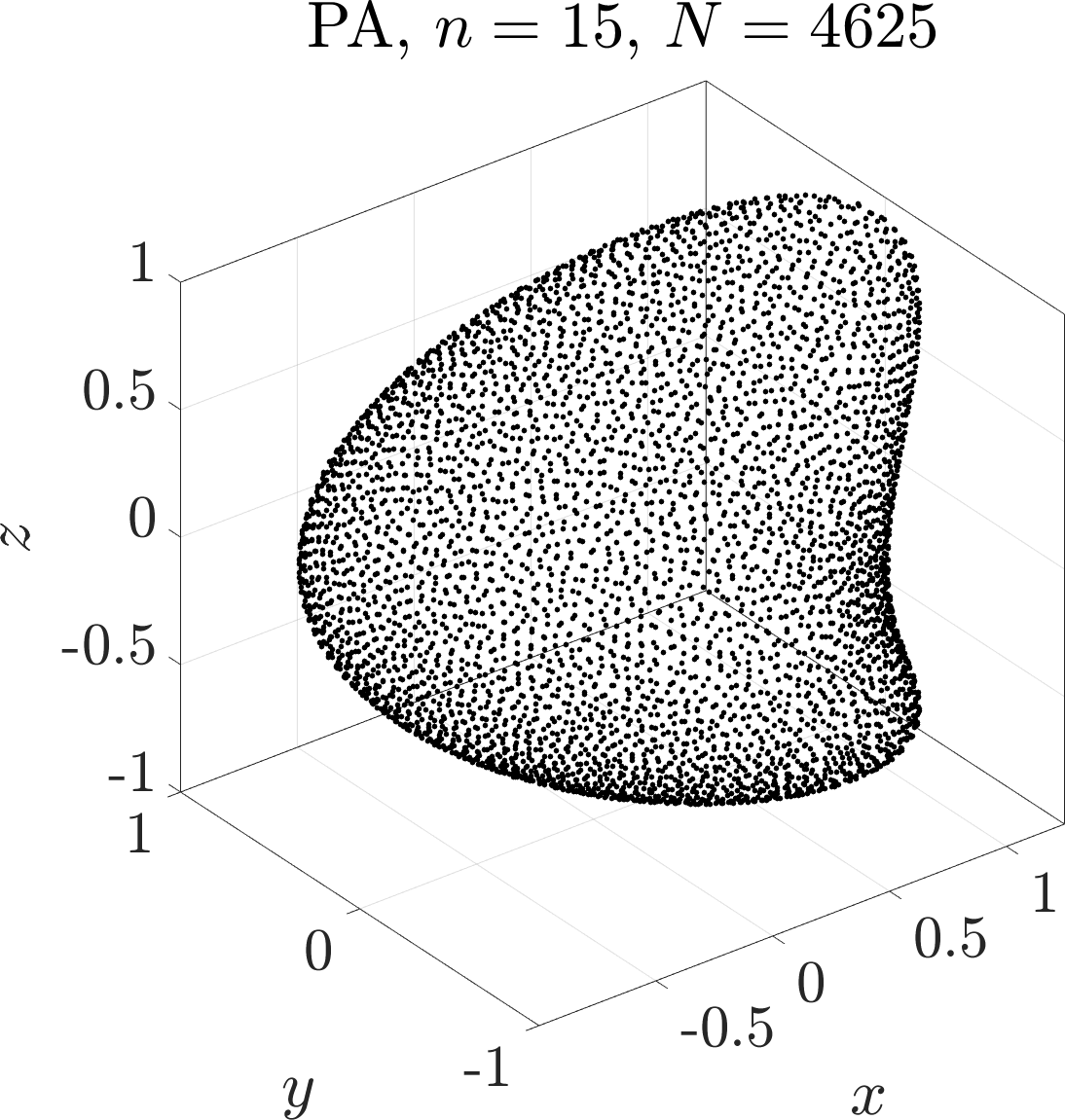}
  \caption{Comparison of different algorithms on a 3D heart-like surface from~\eqref{eq:heart-surf} sampled with $h = 0.05$}
  \label{fig:nodes-3d}
\end{figure}

In~\cref{fig:nodes-3d} we can see their performance in 3D on the heart-like surface $\r_h$.
Naive algorithm's performance is similarly poor as in 2D case. Nodes
generated by the supersampling algorithm have bigger gaps only around $(1, 0 ,1)$ and $(1, 0,
-1)$, where partial derivatives of $\r_h$ diverge. In the same way as in the 2D case,
increasing $\tau$ gives better results, but the increase must be dependent on $h$ to achieve
good quality in all cases. The proposed algorithm also performs worse near
those points. However, it can handle such problematic areas automatically and
provides stable results regardless the $h$. A more in-depth analysis of
minimal and maximal spacing for various $h$ is presented in~\cref{sec:min-max-spacing}.

To further analyze local regularity, we examine the distribution of distances to
nearest neighbors. For each node $\p_i$ we find $c$
nearest neighbors $\p_{i, j}, j = 1, 2, \dots c$ and
calculate the average distance between them $\bar{d}_i = \frac{1}{c} \sum_{j =
1}^{c}\|\p_i - \p_{i, j}\|$. We also calculate the maximum and minimum distances
for each point, $d_i^{\text{min}}$ and  $d_i^{\text{max}}$, where
\begin{equation}
d_i^{\text{min}} = \min_{j=1, \dots c} \|\p_i - \p_{i, j}\|, \qquad
d_i^{\text{max}} = \max_{j=1, \dots c} \|\p_i - \p_{i, j}\|.
\end{equation}
with $c = 2$ for $d_\Xi = 1$ and $c = 3$ for $d_\Xi = 2$.
The results are presented in terms of normalized distances ($\bar{d}'$,
$(d_i^{\text{min}})'$, etc.) which are scaled by $h$, e.g.\
\begin{equation}
  \bar{d}' = \bar{d}/h.
\end{equation}
Numerical results of this analysis are shown in~\cref{tab:node-distr}.

\begin{figure}[h]
  \centering
  \includegraphics[width=0.32\linewidth]{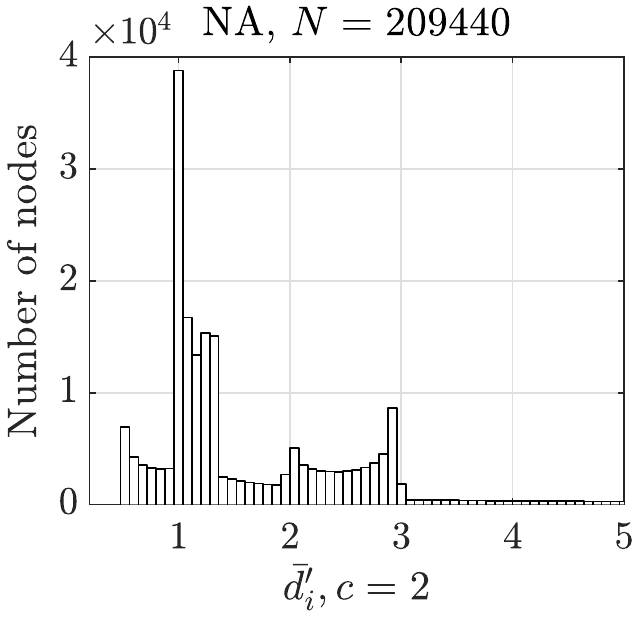}
  \includegraphics[width=0.32\linewidth]{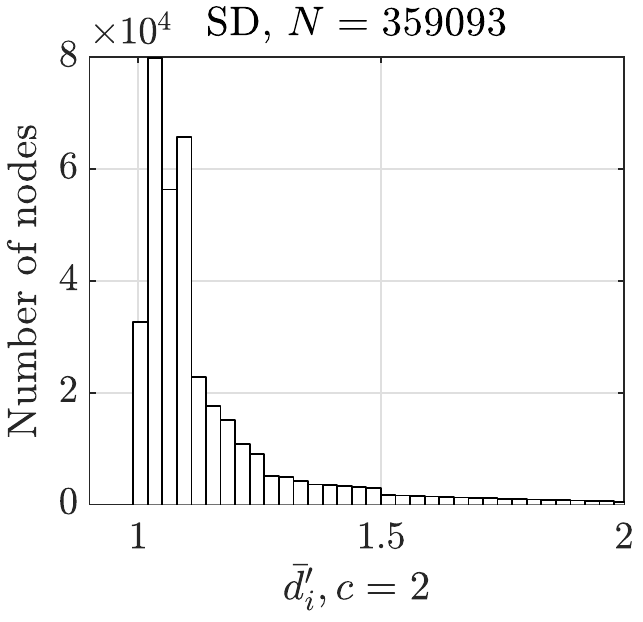}
  \includegraphics[width=0.34\linewidth]{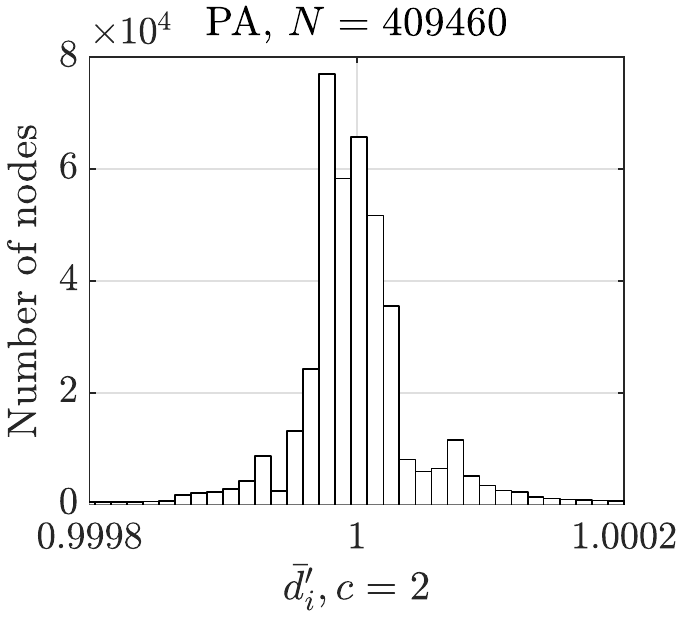}
  \caption{Histogram of normalized average distances to 2 nearest neighbors a polar curve
  from~\eqref{eq:polar-curve} sampled with $h = 0.00003$. $\tau = 5, n = 2$.
  Note the different scaling on the horizontal axes.}
  \label{fig:neighbor-hist-2d}
\end{figure}

\Cref{fig:neighbor-hist-2d} shows the distribution of normalized average
distances for 2 nearest neighbors on the polar curve. The supersampling
algorithm and the proposed algorithm perform much better than the naive algorithm, with
distribution of nodes generated by the proposed algorithm being of higher quality.
The proposed algorithm produces distribution with normalized mean distance
much closer to the target value
$1$ in comparison to the supersampling algorithm. Furthermore, the standard
deviation of $\bar{d}'$ of nodes generated by the proposed algorithm is a few orders of magnitude smaller
than the standard deviation of $\bar{d}'$ of nodes generated by the supersampling algorithm
(see~\cref{tab:node-distr}). Both algorithms show similar mean
difference between the maximum and minimum distance, with the supersampling
algorithm performing slightly better in this aspect.

\begin{figure}[h]
  \centering
  \includegraphics[width=0.32\linewidth]{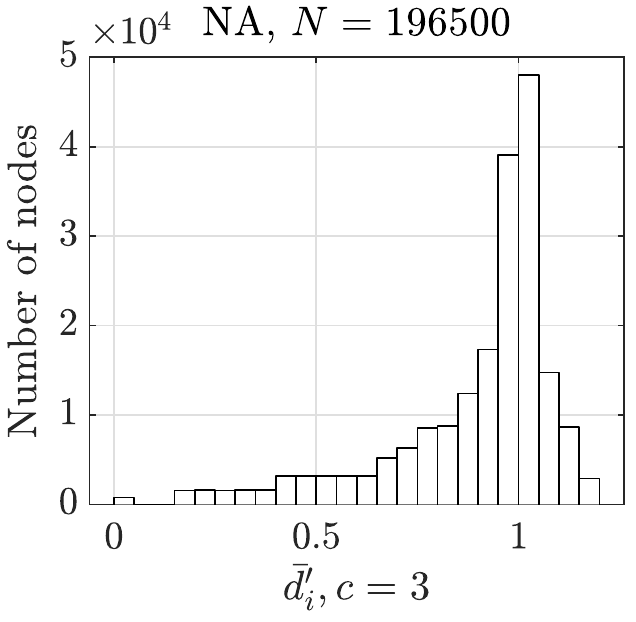}
  \includegraphics[width=0.32\linewidth]{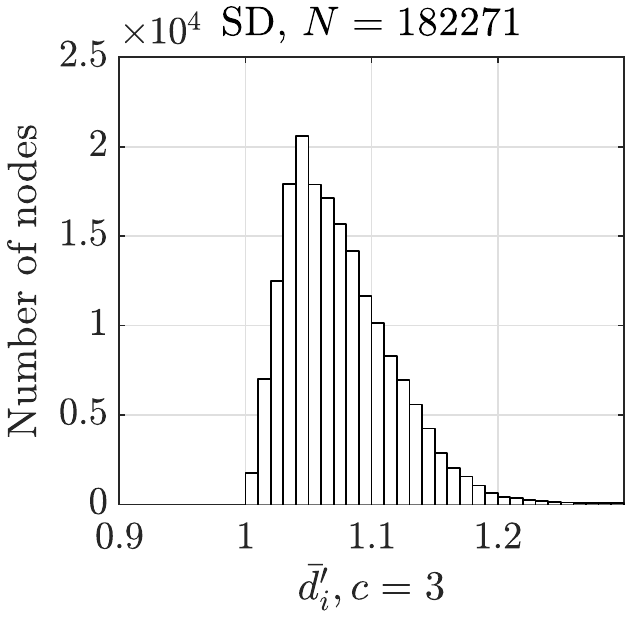}
  \includegraphics[width=0.32\linewidth]{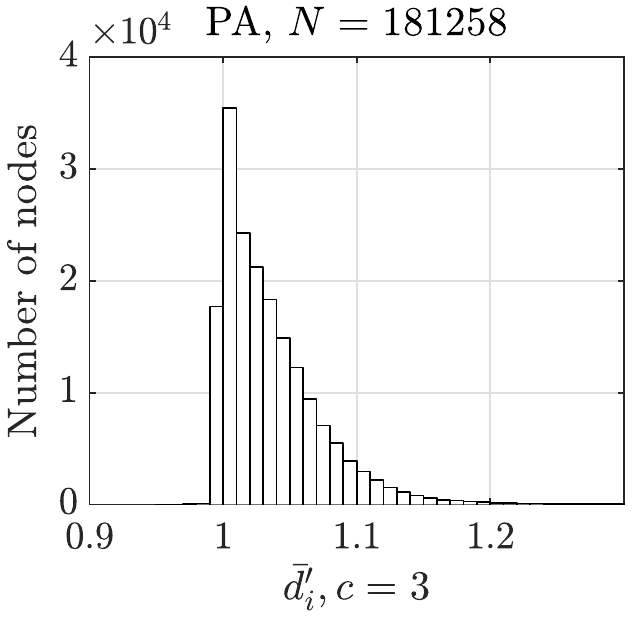}
  \caption{Histogram of normalized average distances to 3 nearest neighbors on a
    heart-like surface from~\eqref{eq:heart-surf} sampled with $h = 0.008$.}
  \label{fig:neighbor-hist-3d}
\end{figure}

\Cref{fig:neighbor-hist-3d} shows the same plots for 3 neighbors on the
heart-like surface. %Zakaj samo 3? Ekvivalentno bi bilo 6, kot v interior 2D fillu? Zato ker so tako rezultati precej boljsi, lahko vecje stevilo, ce bodo recenzenti hoteli
The results
are also similar to the $d_\Xi = 1$ case, the
proposed algorithm has a more favorable mean and standard deviation of  $\bar{d}'$ than supersampling algorithm,
while the naive algorithm is much worse than both of them. However, the differences
between the proposed and supersampling algorithms are much smaller than in the previous
analysis. It is important to
note that increasing the supersampling parameter $\tau$ would, to some degree,
improve the results of supersampling, however at greater computational cost.

\begin{table}[h]
  \renewcommand{\arraystretch}{1.2}
  \centering
  \caption{Numerical quantities related to local regularity. 2D polar curve
  from~\eqref{eq:polar-curve} sampled with $h = 0.00003$ and 3D heart-like
surface from~\eqref{eq:heart-surf} sampled with $h = 0.004$.}
  \label{tab:node-distr}
  \begin{tabular}{c|c|c|c|c}
    case & alg. & $\operatorname{mean}\bar{d}'_i$ & $\operatorname{std}\bar{d}'_i$ &
    $\operatorname{mean}\left(\left(d_i^{\text{max}}\right)' - \left(d_i^{\text{min}}\right)'\right)$ \\ \hline \hline
    \multirow{3}{*}{$d_\Xi = 1$} &
    \rule[-2px]{0px}{14px} PA & 1.0001 & $5.1483 \times 10^{-4}$ & $1.1136
    \times
    10^{-10}$  \\
    & SD & 1.1403 & 0.1715 & $5.7655 \times 10^{-9}$  \\
    & NA & 1.9550 & 1.8386 & $2.7922 \times 10^{-8}$  \\ \hline
    \multirow{3}{*}{$d_\Xi = 2$} &
    \rule[-2px]{0px}{14px} PA  & 1.0357 & 0.0374 & $3.8888 \times 10^{-4}$  \\
    & SD & 1.0764 & 0.0473 & $3.3444 \times 10^{-4}$  \\
    & NA & 0.8946 & 0.2086 & 0.0013  \\
  \end{tabular}
\end{table}

\begin{figure}[h]
  \centering
  \includegraphics[width=0.32\linewidth]{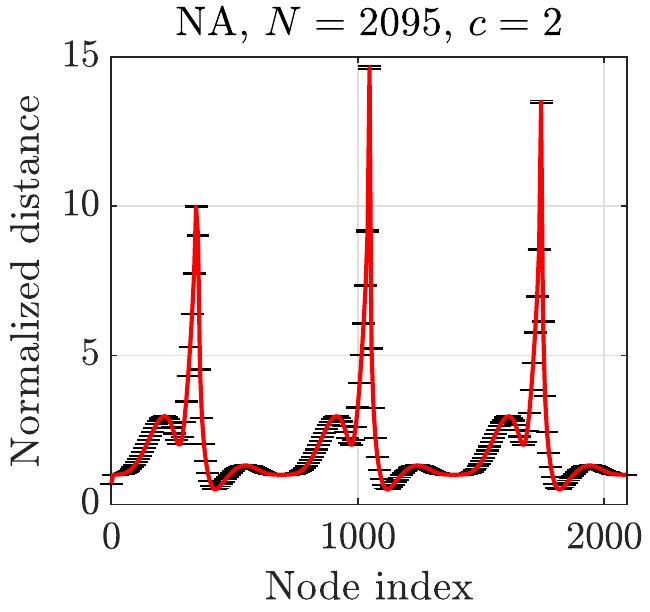}
  \includegraphics[width=0.32\linewidth]{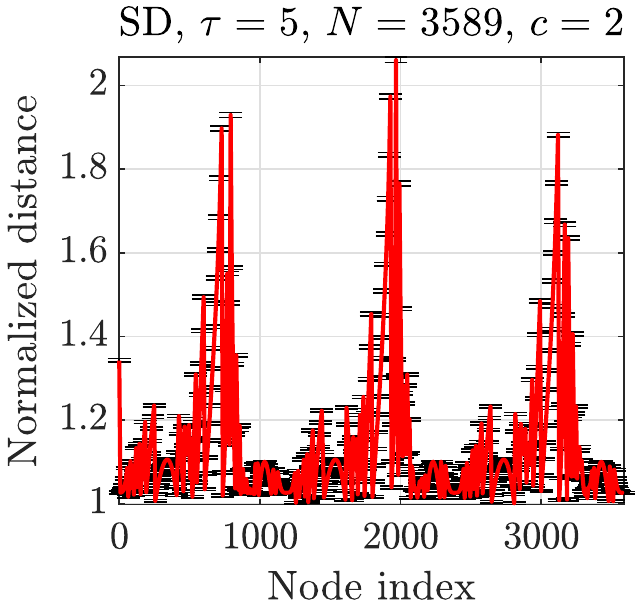}
  \includegraphics[width=0.325\linewidth]{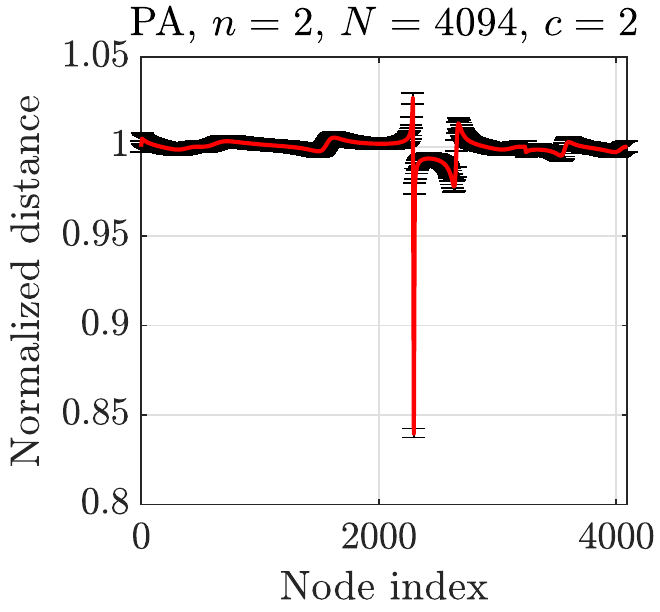}
  \caption{Graph of normalized average distances to 2 nearest neighbors on a 2D
    polar curve from~\eqref{eq:polar-curve} sampled with $h = 0.003$. Error bars
    show minimal and maximal normalized distances to 2 nearest neighbors. Note the different
    $y$-axis range in the plots.}
  \label{fig:neighbor-graph-2d}
\end{figure}

A graph of $\bar{d}_i'$ for each node $\p_i$ with error bars showing
$(d_i^{\text{max}})'$ and $(d_i^{\text{min}})'$
is shown in~\cref{fig:neighbor-graph-2d} for the polar curve.
From this graph we can see that all three algorithms have degraded performance
only on certain points on the curve, where derivatives of $\r_p$ rapidly
increase. Nevertheless, the proposed algorithm handles curves with large
derivatives better.

\begin{figure}[h]
  \centering
  \includegraphics[width=0.32\linewidth]{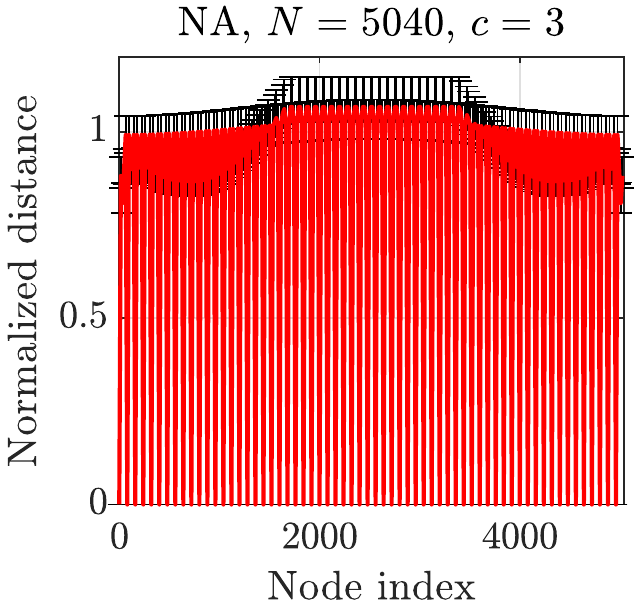}
  \includegraphics[width=0.32\linewidth]{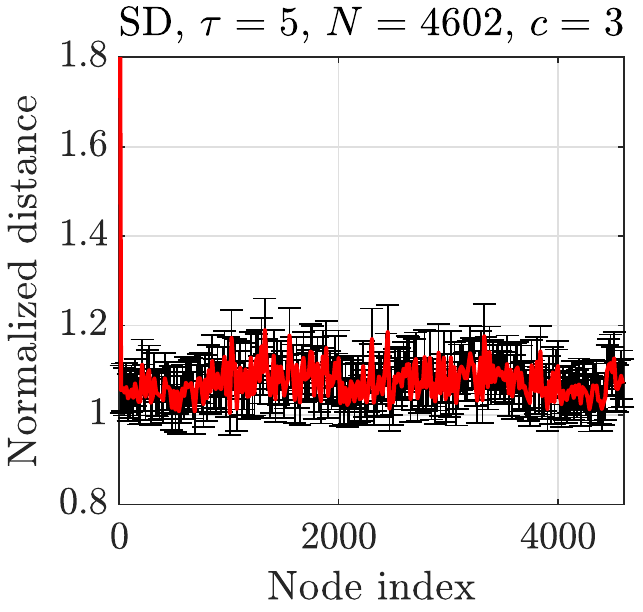}
  \includegraphics[width=0.32\linewidth]{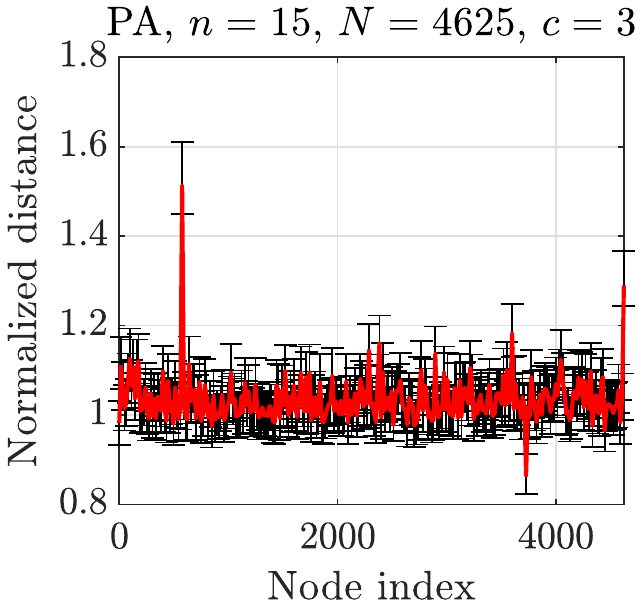}
  \caption{Graph of normalized average distances to 3 nearest neighbors on a 3D
    heart-like surface from~\eqref{eq:heart-surf} sampled with $h = 0.5$. Error
    bars show minimal and maximal normalized distances to 3 nearest neighbors. Note the different
  $y$-axis range in the plots.}
  \label{fig:neighbor-graph-3d}
\end{figure}

\Cref{fig:neighbor-graph-3d} shows the same analysis as~\cref{fig:neighbor-graph-2d}
for $d_\Xi = 2$ case. The differences between the supersampling and the proposed
algorithm decrease and it can be seen that the supersampling algorithm actually has fewer outliers with
high valued $\bar{d}_i'$ and $(d_i^{\text{max}})'$ or low valued $(d_i^{\text{min}})'$.
%All the outliers in both algorithms are however still in the acceptable range.

\subsection{Minimal and maximal spacing requirements}
\label{sec:min-max-spacing}
Minimal and maximal spacing guarantees are in principle inherited from the underlying
algorithm used to discretize $\Xi$, but they are distorted by application of $\r$ in
the naive algorithm. Supersampling algorithm ensures directly that minimal spacing $h$ is
respected in its decimation step. The proposed algorithm uses $\hat{h}_{i,j}$
instead of $h$ to check for distance violations,
and this introduces an error caused by using linear Taylor's
approximation. The exact spacing at point $\p = r(\xib)$ is equal to $h(\p)$, while
the actual computed spacing is equal to
$\hh(\xib, \s) = \|\r(\xib + (h(\p)/ \| \nabla \r (\xib) \s \|) \s) - \r(\xib)\|$.
We wish to estimate the error $\Delta h(\xib, \s) = h(\p) - \hh(\xib, \s)$,
specifically, we would like upper bounds of the form $|\Delta h| \le M$.
Three types of bounds are of interest:
where $M$ depends on $\xib$ and $\s$, where $M$ depends only on $\xi$ and where $M$
is independent and the bound is global.

\begin{proposition}
  \label{prop:bounds}
  The following estimates hold for the error of local node spacing radius
  due to linear approximation in~\eqref{eq:tay-lin}:
  \begin{align}
  |\Delta h(\xib, \s)| &\le \frac{\sqrt{d_\Xi}}{2} h(\p)^2 \frac{ \displaystyle
    \max_{i=1,\ldots,d_\Xi} \max_{\theta \in [0, \alpha]}
    \left|\s^\T (\nabla\nabla r_i)(\xib+\theta \s)\s\right|}{\|\nabla \r(\xib) \s\|^2},
    \alpha = \frac{h(\p)}{ \| \nabla \r (\xib) \s \|}, \\
  |\Delta h(\xib)| &\le \frac{\sqrt{d_\Xi}}{2} h(\p)^2 \frac{ \displaystyle
    \max_{i=1,\ldots,d_\Xi} \max_{\b \zeta \in \bar B(\xib, \rho_{\xib})}
    \sigma_1((\nabla\nabla r_i)(\b \zeta))}{\sigma_{d_\Xi}(\nabla \r(\xib))^2}, \;
  \rho_{\xib} = \frac{h(\p)}{\sigma_{d_\Xi}(\nabla \r(\xib))},  \\
  |\Delta h| &\le \frac{\sqrt{d_\Xi}}{2} h_M^2 \frac{\sigma_{1,M} (\nabla\nabla\r)}{
    \sigma_{d_\Xi, m}^2(\nabla\r)},
  \end{align}
  where \begin{align}
  h_M &= \max_{\xib \in \Xi} h(\p),\\
  \sigma_{1,M} (\nabla\nabla\r) &= \max_{i=1,\ldots,d_\Xi} \max_{\xib \in \Xi}
  \sigma_1((\nabla\nabla r_i)(\xib)), \\
  \sigma_{d_\Xi, m}(\nabla\r) &= \min_{\xib \in \Xi} \sigma_{d_\Xi}(\nabla \r(\xib)),
  \end{align}
  and $\sigma_i(A)$ denotes the $i$-th largest singular value of $A$.

  In particular, this means that the relative error in spacing
  $|\Delta h|/h$ decreases linearly with $h$ for well behaved $\r$ and the algorithm for
  placing points on surfaces asymptotically retains the minimal spacing and quasi-uniformity bounds
  of the underlying algorithm for flat space.
\end{proposition}
%\begin{remark}
%  The bounds involving singular values can be further simplified to a more computationally
%  friendly form. Notably, if all second partial derivatives of $\r$ are
%  bounded as
%  \begin{equation}
%  \max_{\xib \in \Xi}\left|\dpar{^2 r_i}{\xi_j \partial \xi_k}(\xib)\right| \le M,
%  \end{equation}
%  we can use the fact that $\sigma_i(A) \le \|A\|_F$ to obtain
%  \begin{equation}
%  \sigma_{1,M}(\nabla\nabla r) \leq d_{\Xi} M.
%  \end{equation}
%  In a similar fashion, if we know that
%  \begin{equation}
%  \left|\dpar{r_i}{\xi_j}(\xi)\right| \ge m,
%  \end{equation}
%  for some $i, j$, we can transform the expression as $\sigma_{d_\Xi}(A) = \sigma_1(A) / \kappa(A)$
%  and use the bound $\sigma_1(A) \ge a_{ij}$, to obtain
%  \begin{equation}
%  \sigma_{d_\Xi, m}(\nabla\r(\xib)) \geq \frac{m}{\kappa(\nabla \r(\xib))},
%  \end{equation}
%  where $\kappa$ denotes the condition number.
%  This gives us a perhaps more meaningful estimate
%  \begin{equation}
%    WRONG: there is \sigma^2 in the denominator
%  |\Delta h| \le \frac{d_\Xi^{\frac{3}{2}}}{2} h_M^2 \frac{M}{m} \kappa(\nabla \r), \quad
%  \kappa(\nabla \r) = \max_{\xib \in \Xi} \kappa(\nabla \r(\xib)).
%  \end{equation}
%\end{remark}

\Cref{prop:bounds} tells us that the relative error due to linear approximation when computing
$\alpha$ is of order $h$, where the proportionality constant depends on properties of $\nabla \r$
and $\nabla\nabla \r$. Its proof is given in~\Cref{app:proof}.

To analyze minimal and maximal spacing quantitatively, we can use
standard concepts (see e.g.~\cite{wendland2004scattered,hardin2016comparison}), such as
maximum empty sphere radius and separation distance.
For a set of points $\X = \{x_1, \dots, x_N\} \subseteq \partial
\Omega$, the maximum empty sphere radius is defined as \begin{equation}
  r_{\max, \X} = \sup_{x \in \partial \Omega} \min_{1 \leq j \leq N} \|x - x_j\|
\end{equation}
and the separation distance is defined as \begin{equation}
  r_{\min, \X} = \frac{1}{2} \min_{i \neq j} \| x_i - x_j \|
\end{equation}
For numerical computation we compute $r_{\min, \X}$ exactly using a spatial search
structure, such as a $k$-d tree, and we estimate $r_{\max, \X}$ by discretizing the
surface with a much smaller nodal spacing $h$ and finding the maximum empty sphere
with center in one of the generated nodes.

\begin{figure}[h]
  \centering
  \includegraphics[width=0.80\linewidth]{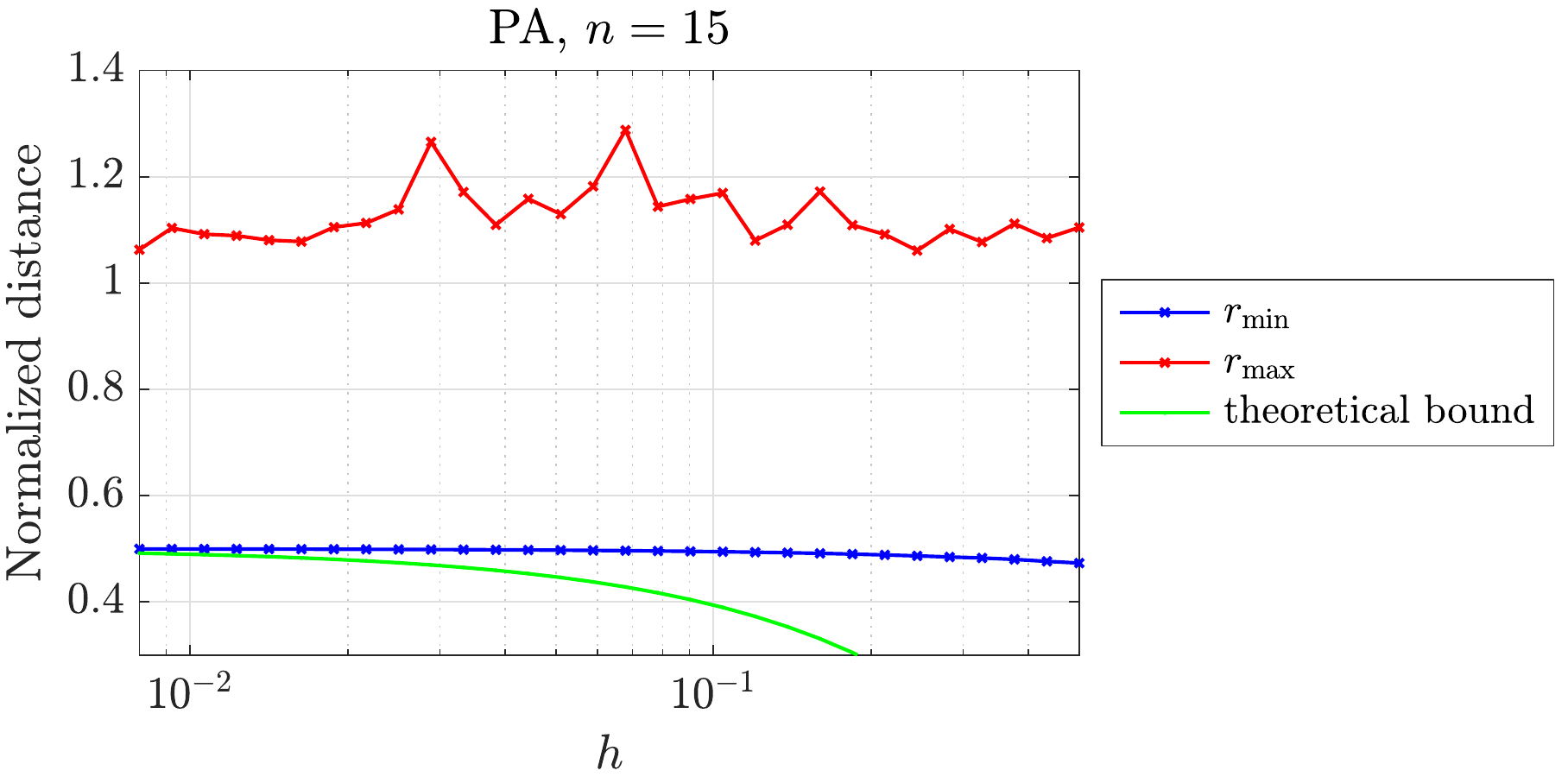}
  \caption{Comparison of practical separation distance and its theoretical lower on nodes generated
  by the proposed algorithm on a torus in 3D.
}
  \label{fig:bound-check-3d}
\end{figure}

According to~\cref{prop:bounds} we conclude that
$r_{\min, \X}/ h \geq 1 + |\Delta h|/ h$ for constant $h$.
This bound has been tested on a torus parameterized with
$\r(\xi_1, \xi_2) = ((\cos\xi_2 + 2)\cos \xi_1, (\cos\xi_2 + 2)\sin \xi_1, \sin \xi_2),
\xi_1, \xi_2 \in [0, 2 \pi]$. By using a computer algebra system, we can calculate that
$\sigma_{1,M} (\nabla\nabla\r) / \sigma_{d_\Xi, m}^2(\nabla\r) = 3$ and compare
this with practical results, shown in~\cref{fig:bound-check-3d}. We can see
that the $r_{\min, \X}$ obeys its lower bound.

\begin{figure}[h]
  \centering
  \includegraphics[width=0.32\linewidth]{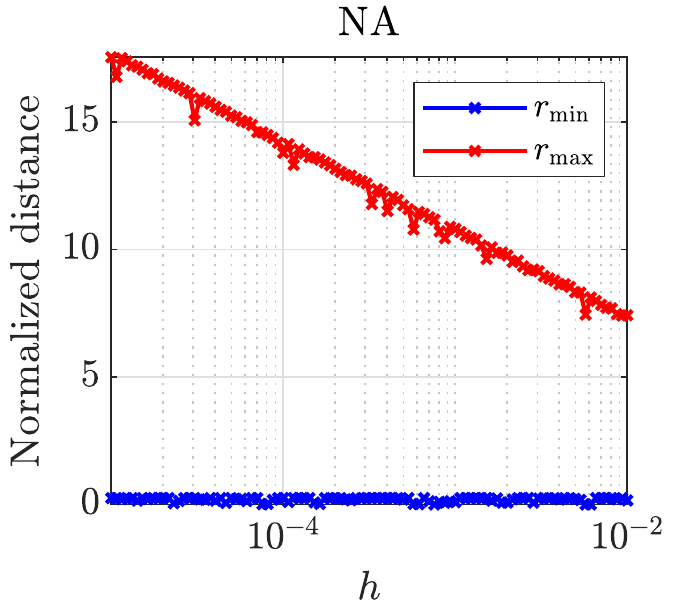}
  \includegraphics[width=0.32\linewidth]{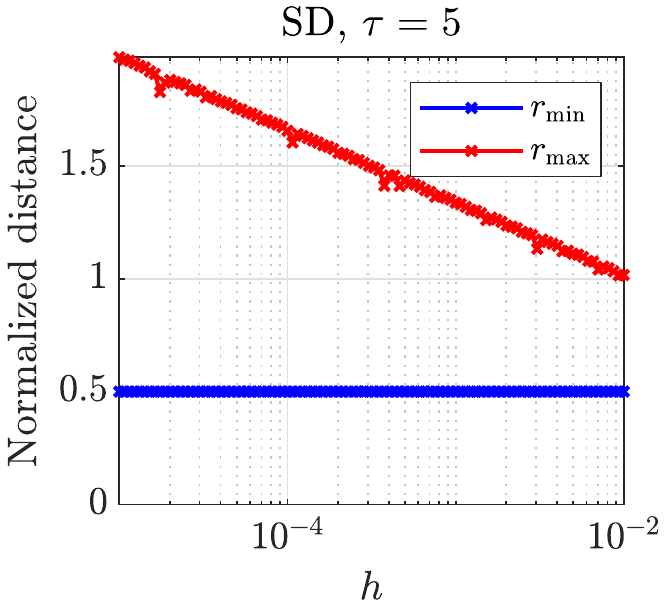}
  \includegraphics[width=0.32\linewidth]{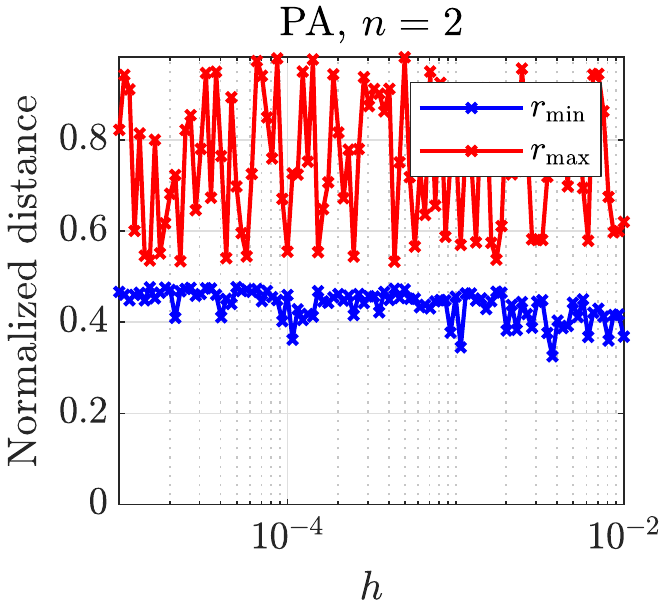}
  \caption{Graph of separation distance and maximum empty sphere radius for different $h$ on the
  polar curve $\r_p$~\eqref{eq:polar-curve}. Note the different $y$-axis spans on the plots. }
  \label{fig:quasi-uni-2d}
\end{figure}

\Cref{fig:quasi-uni-2d} shows a graph of maximum empty sphere radius and
normalized separation distance of the polar curve $\r_p$.
For the naive and supersampling algorithms, normalized maximum empty sphere radius is increasing
with decreasing $h$, because neither of them is able to adapt to the varying value of partial
derivatives, and both perform poorly when $\nabla \r_p$ has high values and when $h$ is smaller.
The proposed algorithm scales much better, i.e.\ $r_{\max}$ remains relatively
stable,
but does not have a strict separation distance minimum of $0.5$, due to the linear
approximation error, as discussed previously.

%, which can be a disadvantage. If a strict
%separation distance is needed, PA can be tweaked to support it.

%Create Chart with Two y-Axes - nima smisla, spodnje meje pri NA in SD so prakticno konstante, samo vec zmede bi bilo

\begin{figure}[h]
  \centering
  \includegraphics[width=0.32\linewidth]{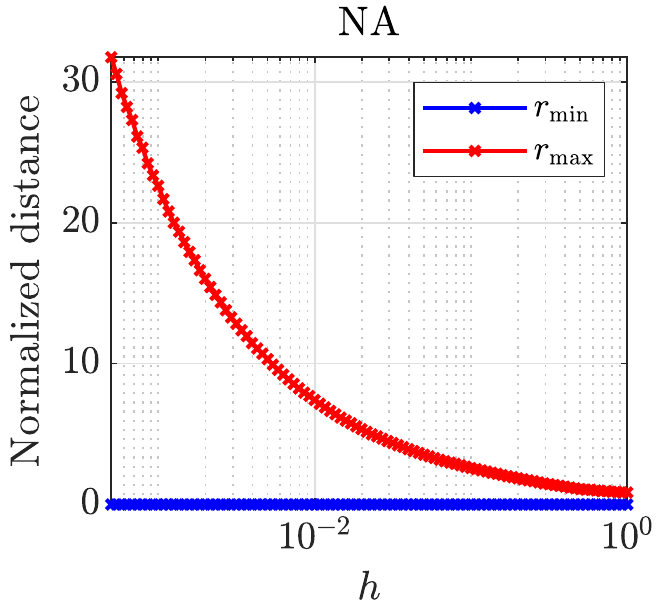}
  \includegraphics[width=0.32\linewidth]{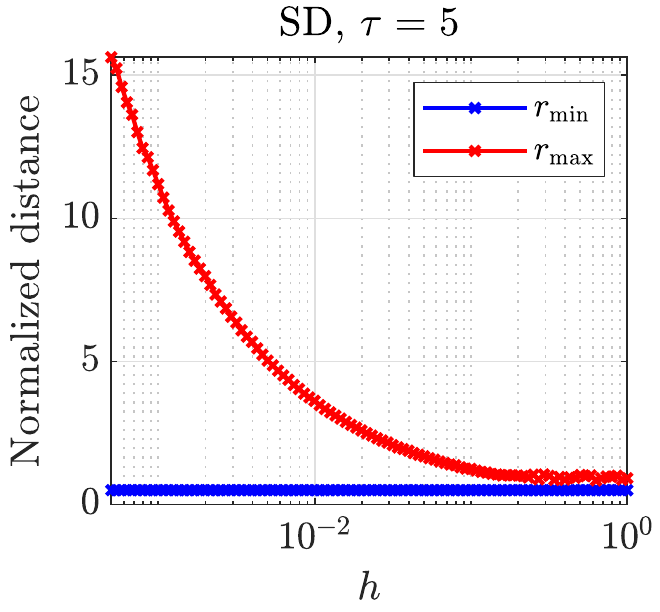}
  \includegraphics[width=0.32\linewidth]{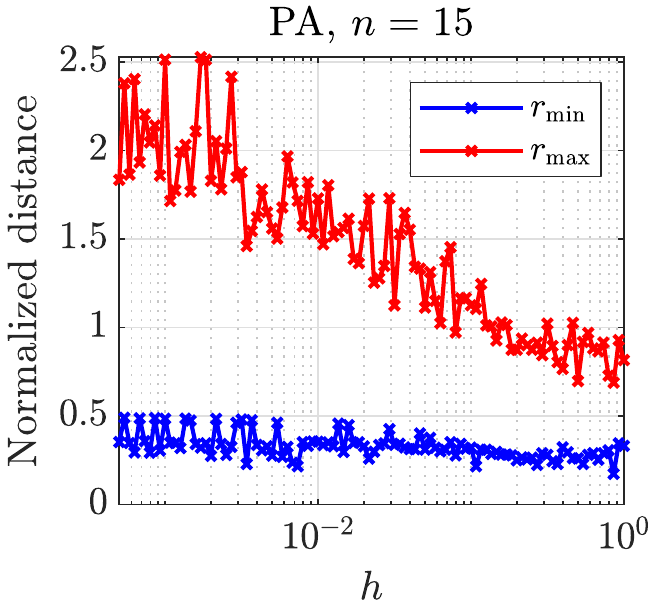}
  \caption{Graph of separation distance and maximum empty sphere radius for different $h$ for the
  heart-like surface $\r_h$~\eqref{eq:heart-surf}. Note the different $y$-axis spans on the plots.}
  \label{fig:quasi-uni-3d}
\end{figure}

In~\cref{fig:quasi-uni-3d} we can see the same graph for the heart-like surface
$\r_h$~\eqref{eq:heart-surf}. The results are similar, however $r_{\max}$
is increasing faster and does not seem bounded for nodes
generated by the proposed algorithm. The reason is similar to the one in 2D, only
in
this case, higher order partial derivatives are increasing faster. If this is
causing problems, the proposed algorithm can be improved by using higher orders in
the Taylor expansion discussed in~\cref{sec:alg}. The performance of the supersampling
algorithm can also be improved by using the higher value of the supersampling
parameter $\tau$, which delays the problem until even smaller $h$. To mitigate this
issue, $\tau$ should be appropriately increased every time $h$ is decreased.

\subsection{Spatial variability}

\begin{figure}[h]
  \centering
  \includegraphics[width=0.49\linewidth]{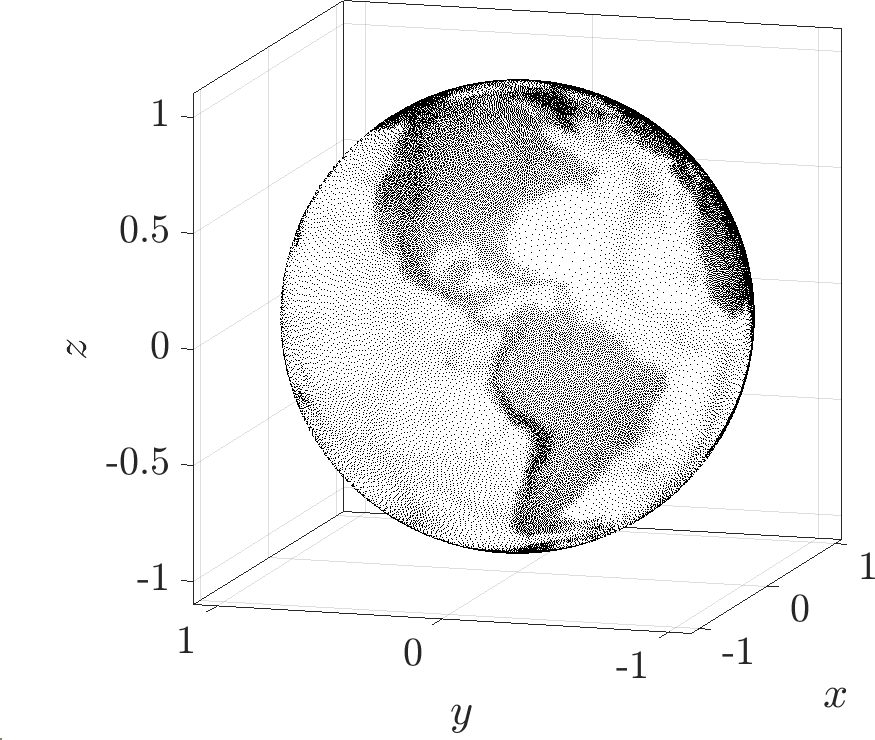}
  \includegraphics[width=0.49\linewidth]{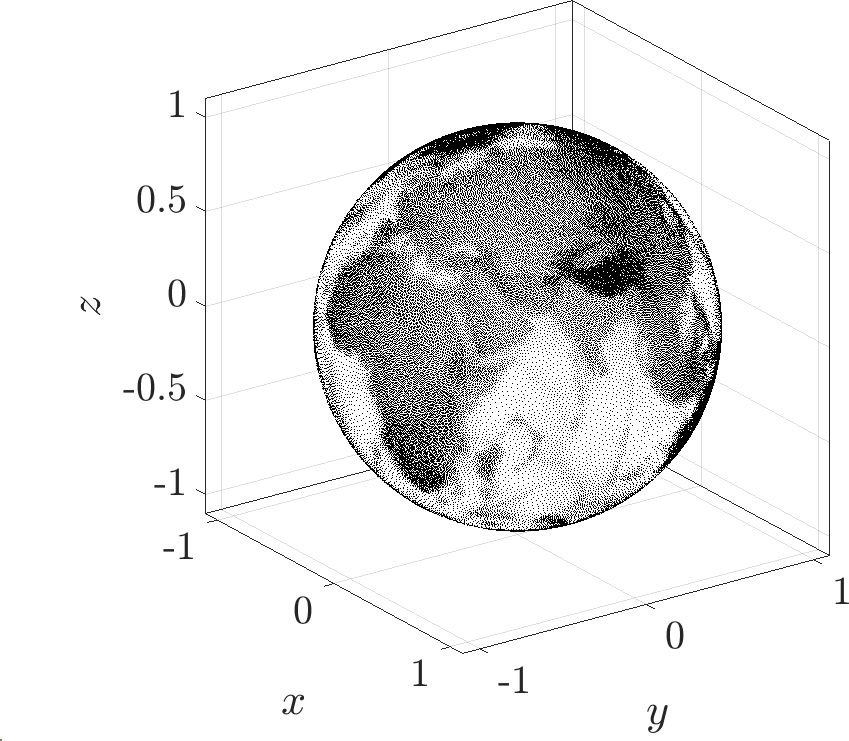}
  \caption{Spherical Earth model sampled proportionally with altitude, $N = 100\,457$.}
  \label{fig:earth}
\end{figure}

An important feature of the proposed algorithm, which other algorithms discussed here do not posses,
is sampling of parametric surfaces with non-constant spacing functions $h$. As a demonstration,
a spherical model of Earth (using the parametrization in spherical coordinates)
was sampled by the proposed algorithm, with density function $h$ representing
altitudes at different points on Earth. Our implementation of the proposed algorithm 
has generated $100\,457$ nodes in less than $\unit[0.3]{s}$ on an 
\texttt{Intel(R) Xeon(R) CPU E5-2620 v3 @ 2.40GHz}. The necessary data was acquired from
Matlab\textsuperscript{\textregistered}'s \texttt{topo.mat} file, also available from
US National Geophysical Data Center~\cite{topo}. The results can be seen in~\cref{fig:earth}.
% Note that this cannot be done with the naive and supersampling algorithms.
% Authors of the supersampling algorithm have
%assumed only constant density functions their article, although a simple tweak
%can make it compatible with non-constant density functions.

\subsection{Computational complexity and execution time}
Theoretical computational complexity of the proposed algorithm was already analyzed
in~\cref{sec:time-theoretical}, and was derived to be
\begin{equation}
  T_{\text{PA}} = O(nN \log N),
\end{equation}
for the general version with a $k$-d tree spatial structure.
The computational complexity of the naive algorithm is clearly
\begin{equation}
T_{\text{NA}} = O(N),
\end{equation}
since the parameters are generated on a grid inside a box-shaped parametric domain $\Xi$
and mapped to the surface.

The computational complexity of the supersampling-decimation algorithm can be written in terms of
the number of generated parameters. If $N_p$ is the number of parameters that are
generated in the parametric domain $\Xi$ with spacing $h_{\Xi}=h/\gamma$,
the time complexity of the algorithm is
$O(N_p \log(N_p))$, as the mapping of the parameters taken $O(N_p)$ time and
closest node queries take $\log(N_p)$ per node using a $k$-d tree data structure.
To compare this result other time complexities, it would need to be expressed in terms of $N$,
the number of finally accepted nodes. It always holds that $N \leq N_p$,
but relating $N$ to $N_p$ in the form of $N_p = O(f(N))$ is not trivial in general
and depends on the parametrization $\r$.  However, by our choice of $\gamma$ it
can be estimated that the SD algorithm generates
$\tau^{d_\Xi} |\obb \partial \Omega|/|\partial \Omega|$ more points than it returns,
where $|\obb \partial \Omega|$ and $|\partial \Omega|$ are surface areas of the oriented bounding
box and the parameterized surface, respectively. Thus we can informally estimate
\begin{equation}
T_{\text{SD}} \approx O\left(\tau^{d_\Xi} \frac{|\obb \partial \Omega|}{|\partial \Omega|} N
\log\left(\tau^{d_\Xi} \frac{|\obb \partial \Omega|}{|\partial \Omega|}N\right)\right).
\end{equation}
To achieve good quality of nodes, the parameter $\tau$
needs to be as large or larger than $\max_{\p \in \partial \Omega} \|\nabla \r(\p)\|$.
Nonetheless, the supersampling-decimation algorithm can be faster than the proposed algorithm
in many real world cases, while the proposed algorithm offers consistent execution time
over a wider range of cases.

\begin{figure}[h]
  \centering
  \includegraphics[width=0.49\linewidth]{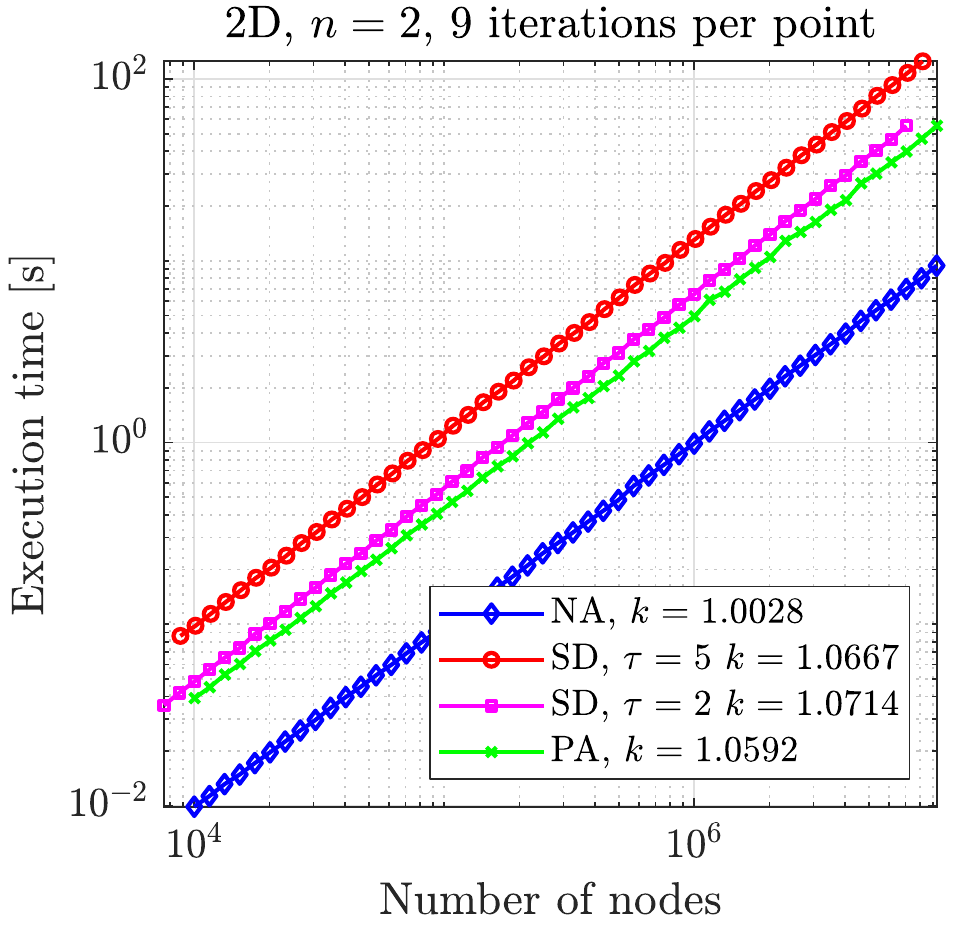}
  \includegraphics[width=0.49\linewidth]{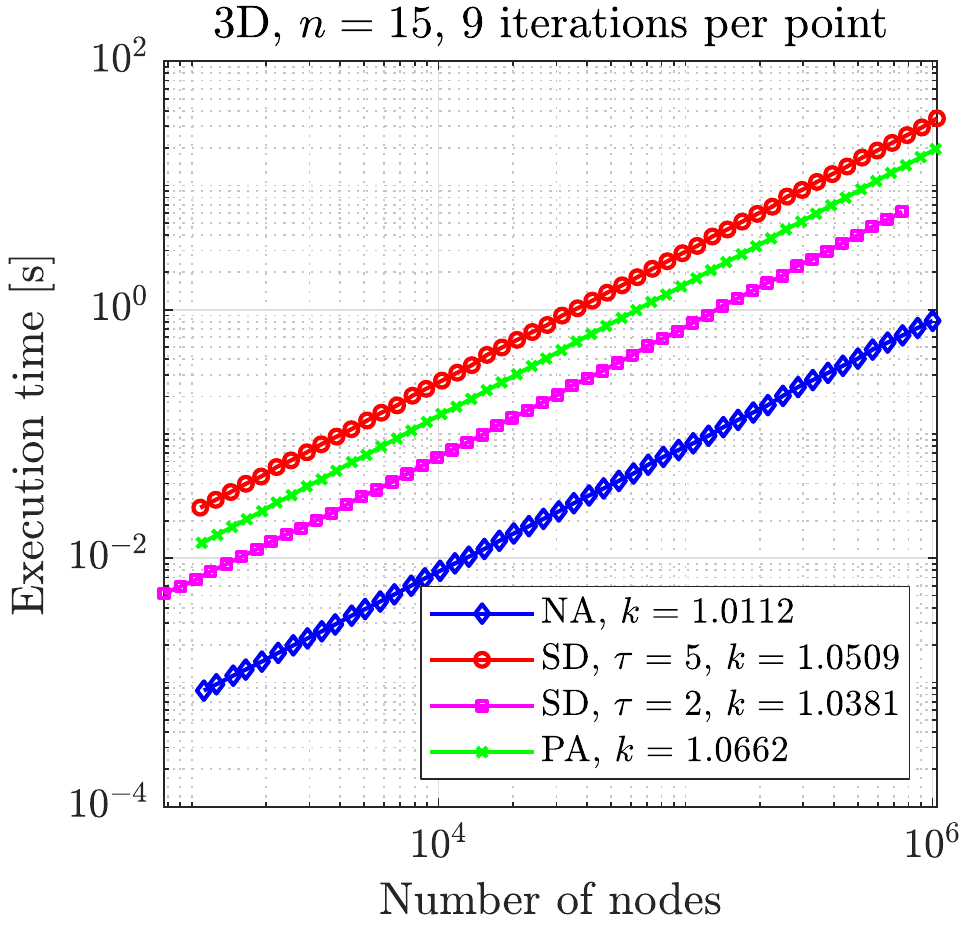}
  \caption{Execution time for different algorithms for different densities in
  2D and 3D. Sampling of 2D polar curve from~\eqref{eq:polar-curve} and 3D
heart-like surface from~\eqref{eq:heart-surf}. The values $k$ represent the estimated line slopes.
}
  %Zaenkrat je tudi v 2D s kd-treejem, lahko naredimo brez, samo verjetno je potem tudi supersampling fer naredit brez, kar je mozno
  \label{fig:exe-time}
\end{figure}

To compare the actual execution time of the three algorithms
they were run with various $h$ on curve $\r_p$~\eqref{eq:polar-curve} and
the heart-like surface $\r_h$~\eqref{eq:heart-surf}.
%All algorithms were implemented in C++ using the Eigen matrix
%library~\cite{eigenweb} and the \texttt{nanoflann} library for $k$-d trees,
%provided by Blanco and Rai~\cite{blanco2014nanoflann}.
The measurements were
done on a machine with an \texttt{Intel(R) Xeon(R) CPU E5-2620 v3 @ 2.40GHz}
processor and 64 GB DDR3 RAM.  Code was compiled with \texttt{g++ (GCC) 8.1.0}
for Linux with \texttt{-std=c++11 -O3 -DNDEBUG} flags.  Measurements for each
data point were executed 9 times and the median time was taken.  The results
are shown in~\cref{fig:exe-time}.

Growth factors obtained from practical results match the theoretical time complexities.
A different choice of $\tau$ can make the supersampling algorithm faster or slower than
proposed algorithm. Both are able to generate $10^6$ nodes in order of a few
seconds in $d_\Xi = 1$ case and in order of few 10s of seconds in $d_\Xi = 2$ case.

\section{Mesh-free numerical analysis example}
\label{sec:num}
In this section we demonstrate the performance of the proposed algorithm in
providing the discretization of domain boundary for meshless solution of PDEs.
Consider a domain $\Omega$ bounded by the polar curve $\r_p$ in 2D
and bounded by the heart-like surface $\r_h$ in 3D.

For the closed form solution we choose $u_2(x, y) = \sin(\pi x) \cos(2 \pi y)$ in 2D and $u_3(x, y,
z) = \sin(\pi x) \cos(2 \pi y) \sin(\frac{1}{2} \pi z)$ in 3D and define
the following Poisson problem in $d$ dimensions (for $d=2, 3$):
\begin{align}
  \hspace{4cm} \nabla^2 u &= f && \text{in $\Omega$},               \hspace{4cm} \\
  \hspace{4cm} u &= u_d && \text{on $\Gamma_e$},                    \hspace{4cm}  \\
  \hspace{4cm} \vec{n} \cdot \nabla u &= g && \text{on $\Gamma_n$}, \hspace{4cm}
\end{align}
where $f$ and $g$ are computed from $u_d$.

The domain $\Omega$ is discretized in two steps. First, boundary
$\partial \Omega$ is discretized using the proposed algorithm. In the second step,
the interior of $\Omega$ is populated with nodes using the algorithm introduced
in~\cite{slak2019generation} where the boundary nodes from the first step are used as
the seed nodes.

Once the domain is fully populated with nodes, a
radial basis function-generated finite differences (RBF-FD) method,
using polyharmonic basis functions augmented with monomials is used.
RBF-FD using polyharmonics augmented with monomials is a promising
mesh-free method that combines the robustness of classical RBF-FD but
circumvents the stagnation errors and achieves high-order accuracy by
leveraging monomial augmentation~\cite{bayona2017augment}.

We used RBF-FD with monomial augmentation up to order $m \in \{2, 4, 6\}$, where we 
chose stencil size $n = 4 \binom{m + 2}{2}$ in 2D and $n = 4 \binom{m + 3}{3}$ in 3D, 
based on guidelines from~\cite{bayona2017augment}, 
to obtain the mesh-free approximations of the differential operators involved.
The resulting sparse system is solved with BiCGSTAB using an ILUT preconditioner.

\begin{figure}[h]
  \centering
  \includegraphics[width=0.32\linewidth]{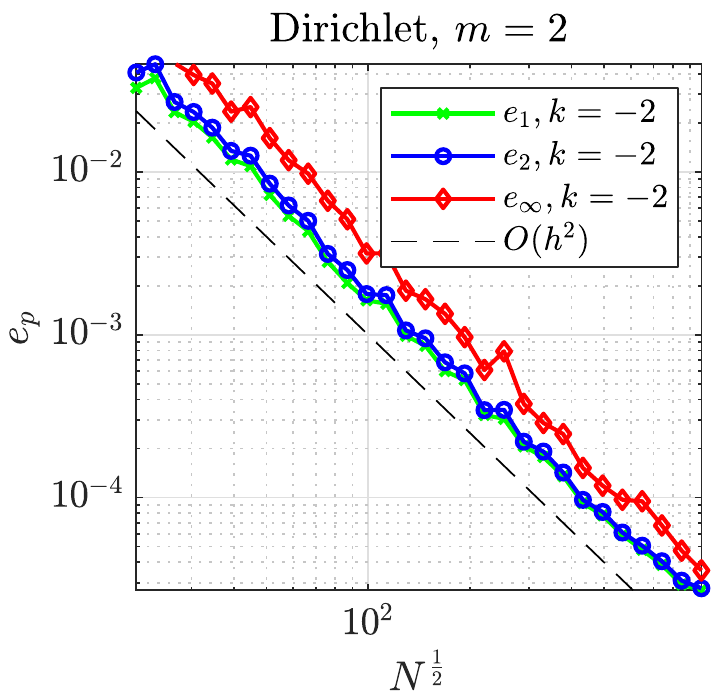}
  \includegraphics[width=0.32\linewidth]{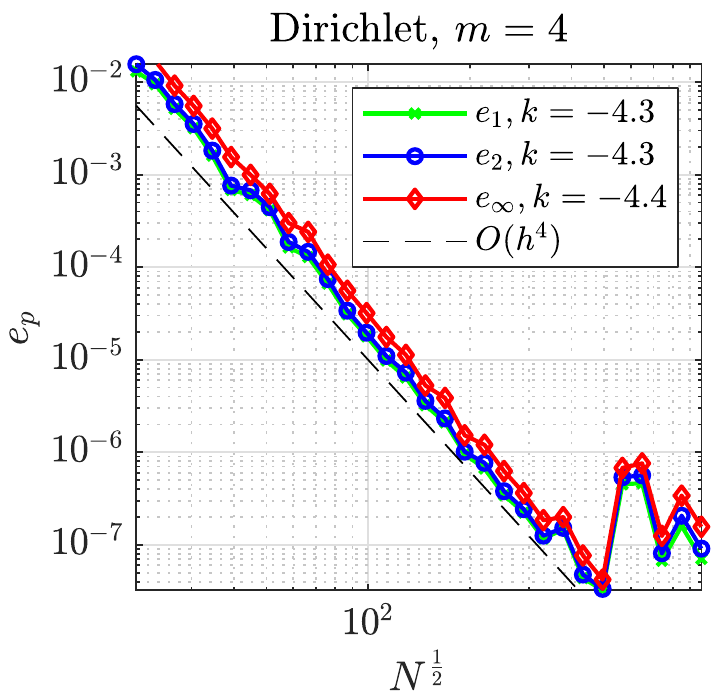}
  \includegraphics[width=0.32\linewidth]{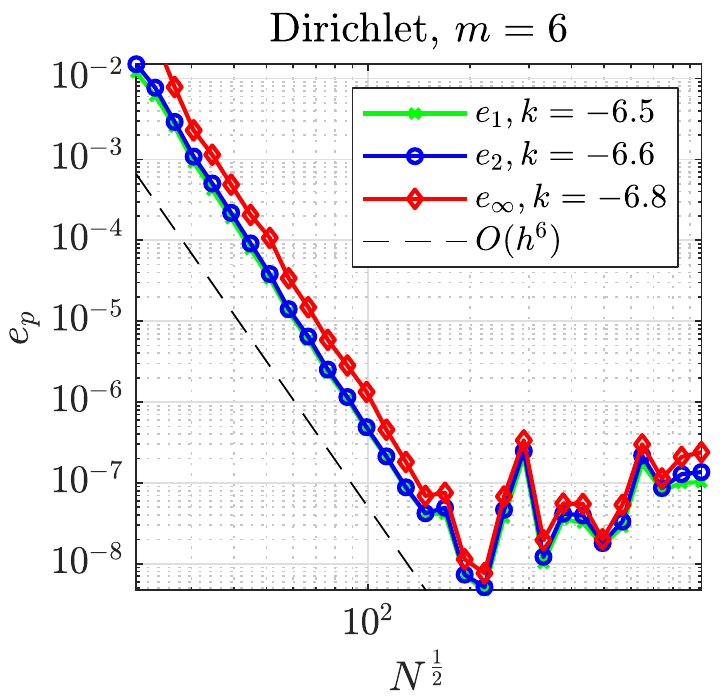}

  \includegraphics[width=0.32\linewidth]{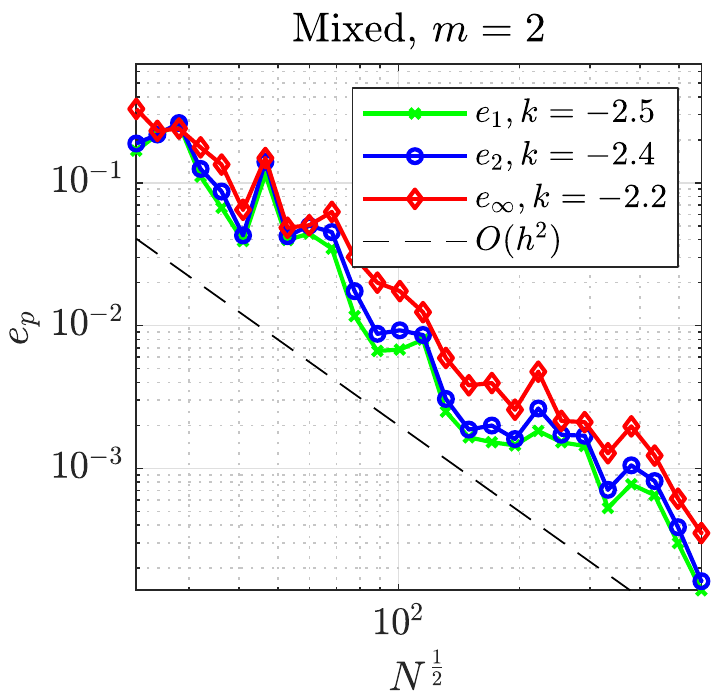}
  \includegraphics[width=0.32\linewidth]{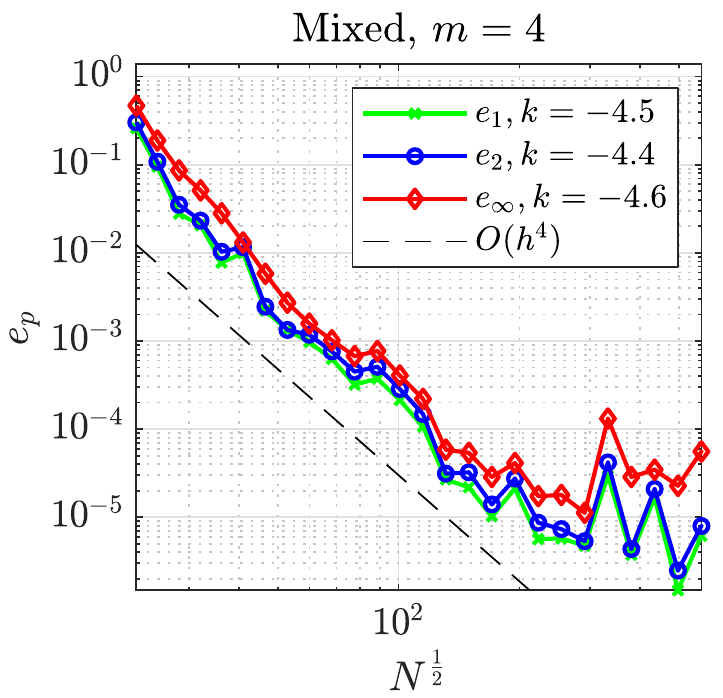}
  \includegraphics[width=0.32\linewidth]{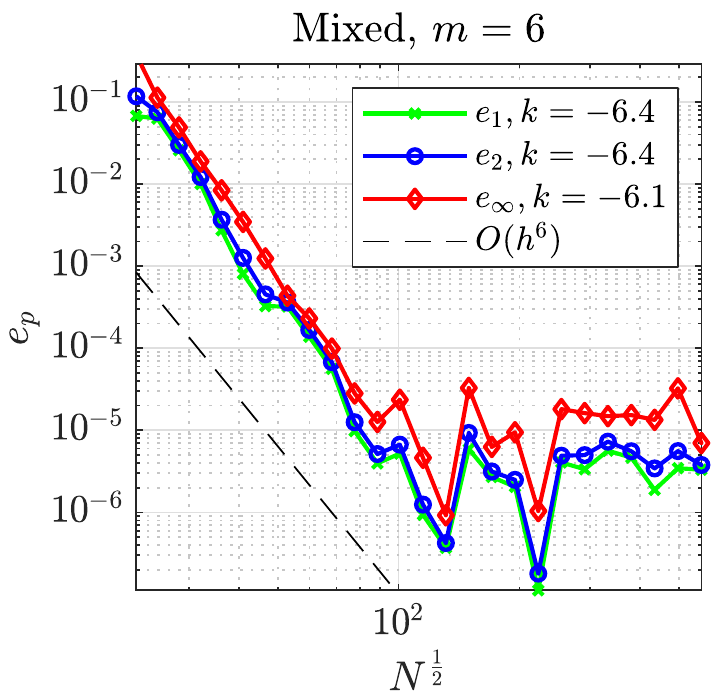}
  \caption{Convergence of Poisson's equation $u(x, y) = \sin(\pi x) \cos(2 \pi
  y)$ with different boundary conditions and monomials up to order $m$ on nodes
  generated by the proposed algorithm and the algorithm presented in~\cite{slak2019generation}
  on 2D polar curve from~\eqref{eq:polar-curve}. The values $k$ represent the estimated line
  slopes until the error starts diverging.}
  \label{fig:poisson-2d}
\end{figure}

\begin{figure}[h]
  \centering
  \includegraphics[width=0.32\linewidth]{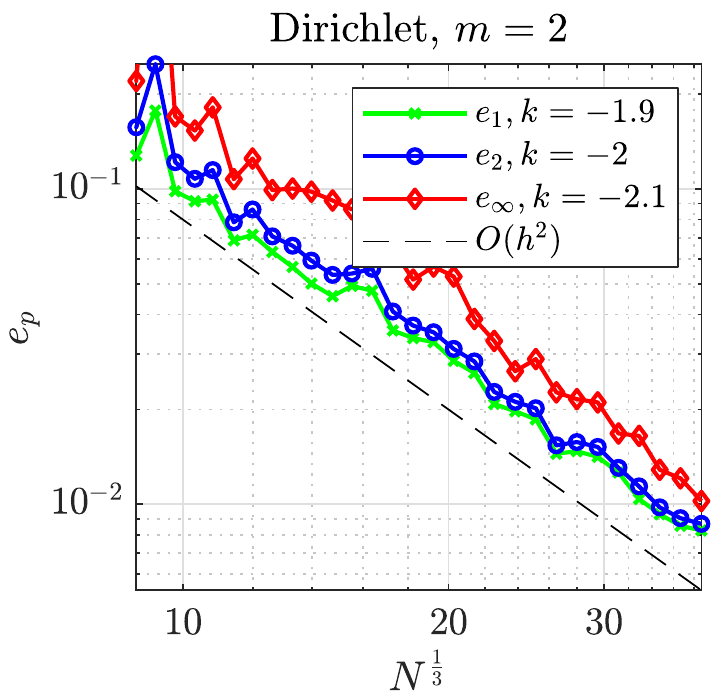}
  \includegraphics[width=0.32\linewidth]{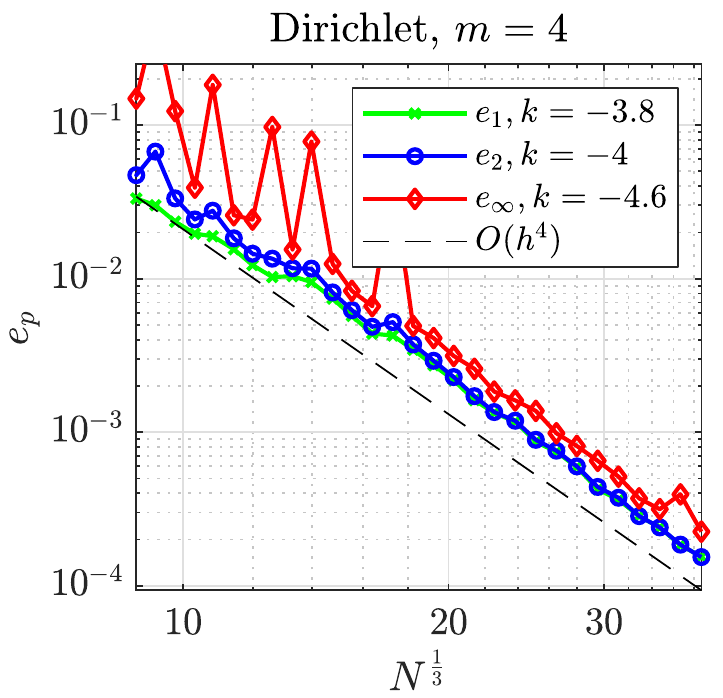}
  \includegraphics[width=0.32\linewidth]{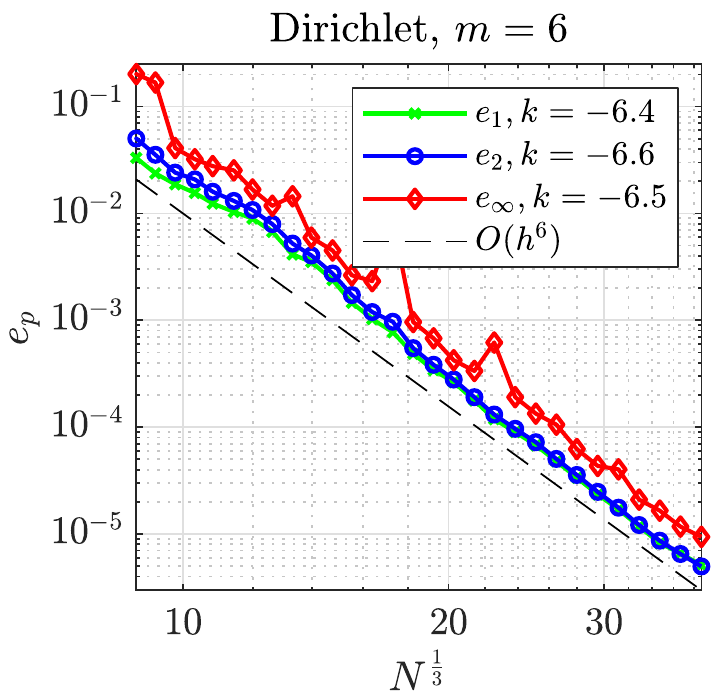}

  \includegraphics[width=0.32\linewidth]{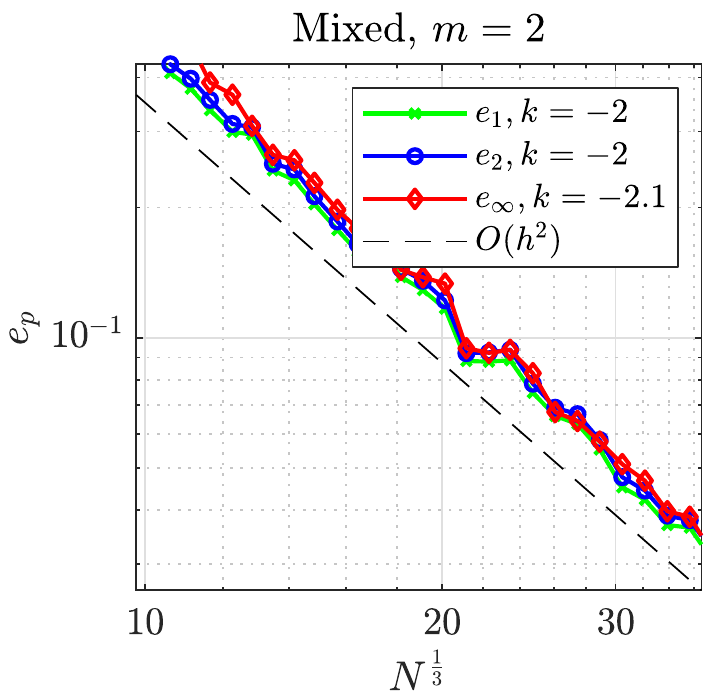}
  \includegraphics[width=0.32\linewidth]{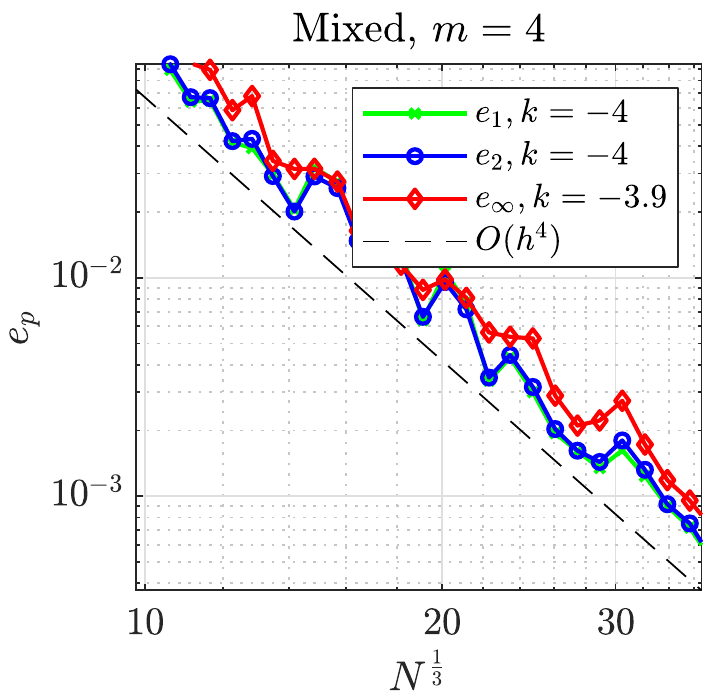}
  \includegraphics[width=0.32\linewidth]{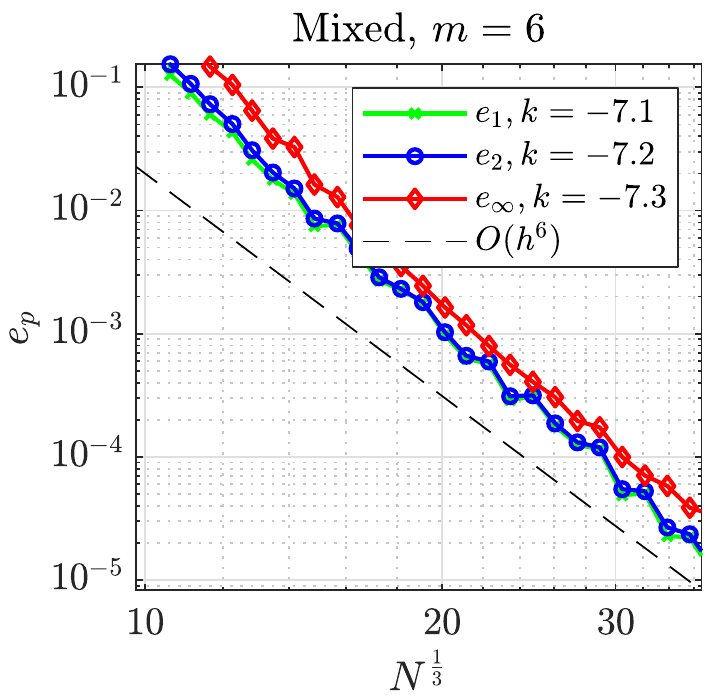}
  \caption{Convergence of Poisson's equation $u(x, y, z) = \sin(\pi x) \cos(2
  \pi y) \sin(\frac{1}{2} \pi z)$ with mixed boundary conditions and
monomials up to order $m$ on nodes generated by the proposed algorithm and the
algorithm presented in~\cite{slak2019generation} on 3D heart-like surface
from~\eqref{eq:heart-surf}. The values $k$
represent the estimated line slopes.}
  \label{fig:poisson-3d}
\end{figure}

Figures~\ref{fig:poisson-2d} and~\ref{fig:poisson-3d} show
the relative discrete $p$-norm error $e_p = \|u - \hat{u}\|_p / \|u\|_p$
for increasing number of nodes in two cases. In one case, only Dirichlet boundary conditions
were used ($\Gamma_e = \partial \Omega$ and $\Gamma_n = \emptyset$) and
in the other case, labeled ``mixed'', Dirichlet boundary conditions were used for boundary nodes
with negative $x$-coordinate and Neumann boundary conditions were used otherwise.

The error behaves as expected in 3D, but starts diverging when the number of
nodes is big in 2D. This is because of finite precision errors discussed by
Flyer et al.~\cite{flyer2016role}. Until that point, the error also behaves as
expected in 2D.

\section{Conclusions}
\label{sec:conc}
A new algorithm for generating nodes on parametric surfaces was developed and
compared with naive parametric sampling and the supersampling-decimation
algorithm~\cite{shankar2018robust}. All three algorithms require that a
parametrization $\r$ of a curve/surface in question is given.

Of the three algorithms, both supersampling and the proposed algorithm produce
quasi-uniform discretizations.
The proposed algorithm requires that $\nabla \r$ is given in addition to $\r$,
which can be problematic, but is usually known in case of closed form parametrization,
RBF models or CAD models.
Contrary to the other two algorithms, the proposed algorithm supports generation of nodes
with variable nodal spacing, and on irregular parametric domains. It
also adapts to parametrizations with variable $\|\nabla \r\|$ without modifications.
We also proved minimal spacing requirements of the proposed algorithm, for both
uniform and variable spacing.

Both the proposed and the supersampling-decimation algorithm support generation of nodes in
arbitrary dimensions. The time complexity of the proposed algorithm is
$O(N \log N)$ to generate $N$ nodes in all cases, while the time complexity of the
supersampling algorithm varies a lot with $\tau$ and the properties of the parametrization~$\r$.
The execution times of the algorithm are comparable, with execution time of
supersampling-decimation algorithm depending heavily on choice of $\tau$.
It takes around \unit[1]{s} to generate $10^6$ nodes in 2D and \unit[10]{s} to generate
$10^6$ nodes in 3D.

Future work is focused on the parallelization of the underlying spatial
generation algorithm, with some successful preliminary results~\cite{parnum} and
analyzing the behavior of the algorithm for surfaces composed of multiple
patches, with the ultimate goal to support automatic meshless discretization of
CAD models.

\section*{Acknowledgments}
The authors would like to acknowledge the financial support of the ARRS
research core funding No.\ P2-0095 and the Young
Researcher program PR-08346.

\appendix
\section{Node spacing error proposition}
\label{app:proof}

\newsiamthmrep{rep_proposition}{Proposition \ref{prop:bounds}}

\begin{rep_proposition}
	The following estimates hold for the error of local node spacing radius
	due to linear approximation in~\eqref{eq:tay-lin}:
	\begin{align}
	|\Delta h(\xib, \s)| &\le \frac{\sqrt{d_\Xi}}{2} h(\p)^2 \frac{ \displaystyle
		\max_{i=1,\ldots,d_\Xi} \max_{\theta \in [0, \alpha]}
		\left|\s^\T (\nabla\nabla r_i)(\xib+\theta \s)\s\right|}{\|\nabla \r(\xib) \s\|^2},
	\alpha = \frac{h(\p)}{ \| \nabla \r (\xib) \s \|}, \\
	|\Delta h(\xib)| &\le \frac{\sqrt{d_\Xi}}{2} h(\p)^2 \frac{ \displaystyle
		\max_{i=1,\ldots,d_\Xi} \max_{\b \zeta \in \bar B(\xib, \rho_{\xib})}
		\sigma_1((\nabla\nabla r_i)(\b \zeta))}{\sigma_{d_\Xi}(\nabla \r(\xib))^2}, \;
	\rho_{\xib} = \frac{h(\p)}{\sigma_{d_\Xi}(\nabla \r(\xib))},  \\
	|\Delta h| &\le \frac{\sqrt{d_\Xi}}{2} h_M^2 \frac{\sigma_{1,M} (\nabla\nabla\r)}{
		\sigma_{d_\Xi, m}^2(\nabla\r)},
	\end{align}
	where \begin{align}
	h_M &= \max_{\xib \in \Xi} h(\p),\\
	\sigma_{1,M} (\nabla\nabla\r) &= \max_{i=1,\ldots,d_\Xi} \max_{\xib \in \Xi}
	\sigma_1((\nabla\nabla r_i)(\xib)), \\
	\sigma_{d_\Xi, m}(\nabla\r) &= \min_{\xib \in \Xi} \sigma_{d_\Xi}(\nabla \r(\xib)),
	\end{align}
	and $\sigma_i(A)$ denotes the $i$-th largest singular value of $A$.

	In particular, this means that the relative error in spacing
	$|\Delta h|/h$ decreases linearly with $h$ for well behaved $\r$ and the algorithm for
	placing points on surfaces asymptotically retains the minimal spacing and quasi-uniformity bounds
	of the underlying algorithm for flat space.
\end{rep_proposition}
%\begin{remark}
%  The bounds involving singular values can be further simplified to a more computationally
%  friendly form. Notably, if all second partial derivatives of $\r$ are
%  bounded as
%  \begin{equation}
%  \max_{\xib \in \Xi}\left|\dpar{^2 r_i}{\xi_j \partial \xi_k}(\xib)\right| \le M,
%  \end{equation}
%  we can use the fact that $\sigma_i(A) \le \|A\|_F$ to obtain
%  \begin{equation}
%  \sigma_{1,M}(\nabla\nabla r) \leq d_{\Xi} M.
%  \end{equation}
%  In a similar fashion, if we know that
%  \begin{equation}
%  \left|\dpar{r_i}{\xi_j}(\xi)\right| \ge m,
%  \end{equation}
%  for some $i, j$, we can transform the expression as $\sigma_{d_\Xi}(A) = \sigma_1(A) / \kappa(A)$
%  and use the bound $\sigma_1(A) \ge a_{ij}$, to obtain
%  \begin{equation}
%  \sigma_{d_\Xi, m}(\nabla\r(\xib)) \geq \frac{m}{\kappa(\nabla \r(\xib))},
%  \end{equation}
%  where $\kappa$ denotes the condition number.
%  This gives us a perhaps more meaningful estimate
%  \begin{equation}
%    WRONG: there is \sigma^2 in the denominator
%  |\Delta h| \le \frac{d_\Xi^{\frac{3}{2}}}{2} h_M^2 \frac{M}{m} \kappa(\nabla \r), \quad
%  \kappa(\nabla \r) = \max_{\xib \in \Xi} \kappa(\nabla \r(\xib)).
%  \end{equation}
%\end{remark}
\begin{proof}
	The following estimate for $a \in \R$ and $b, c \in \R^n$ will be used:
	\begin{equation} \label{eq:dbl-abs}
	|a - \|b+c\|| \le |a - \|b\|| + \|c\|.
	\end{equation}

	%\textbf{Proof:}
	%The following triangle inequality-based bound will be used: $\|b\|-\|c\| \leq \|b+c\| \leq
	%\|b\|+\|c\|$.
	%We differentiate two cases
	%\begin{itemize}
	%  \item $a \leq \|b+c\|$:
	%  \begin{equation}
	%    |a - \|b+c\||  = \|b+c\|-a \leq \|b\|+\|c\| - a = \|b\|-a + \|c\| \leq |\|b\|-a| + \|c\|
	%  \end{equation}
	%  \item $a \geq \|b+c\|$, also implies $a \geq \|b\|$:
	%  \begin{equation}
	%  |a - \|b+c\||  = a-\|b+c\| \leq  a - (\|b\|-\|c\|) \leq a - \|b\| + \|c\| = |a-\|b\|| + \|c\|
	%  \end{equation}
	%\end{itemize}

	The dependence of $\alpha$ on $\xib$ and $\s$ will also be written explicitly, $\alpha(\xib, \s)$.
	We begin to estimate the error locally
	\begin{align}
	|\Delta h(\xib, \s)|  &\le |h(\p) - \hh(\xib, \s)| \\
	&= |h(\p) - \|\r(\etab) - \r(\xib)\||
	= |h(\p) - \|\alpha(\xib, \s) \nabla \r (\xib) \s + \Rem(\xib, \vec s) \| | \\
	&\le |h(\p) - \|\alpha(\xib, \s) \nabla \r (\xib) \s\|| + \|\Rem(\xib, \vec s) \| \\
	&= \left|h(\p) - \frac{h(\p)}{\|\nabla \r (\xib) \s \|}\|\nabla \r (\xib) \s\|\right| +
	\|\Rem(\xib, \vec s) \| \\
	&= \|\Rem(\xib, \vec s) \|,
	\end{align}
	where $\Rem(\xib, \vec s)$ is the remainder of the Taylor approximation
	\begin{equation}
	\Rem(\xib, \vec s) = \r(\etab) - \r(\xib) - \alpha \nabla \r (\xib) \s,
	\end{equation}
	which can also we viewed as an Taylor expansion of $\r(\p + \alpha \s)$ around $\alpha = 0$.
	The remainder can be further estimated component-wise:
	\begin{equation}
	\|\Rem(\xib, \vec s)\| \leq \sqrt{d_{\Xi}} \max_{i=1,\ldots,d_\Xi} \|R_i(\xib, \vec s)\|.
	\end{equation}

	Using the Lagrange form of remainder for each component of $R(\xib)$, we arrive at
	\begin{equation}
	R_i(\xib, \vec s) = \frac{1}{2} h(\p)^2
	\frac{\s^\T (\nabla\nabla r_i)(\xib+\theta \s)\s}{\|\nabla \r (\xib) \s\|^2}
	\end{equation}
	for some $\theta \in [0, \alpha]$, where $\nabla\nabla r_i$ is the Hessian matrix of
	$r_i$. Thus we can bound each component as
	\begin{equation} \label{eq:component-bound}
	|R_i(\xib, \vec s)| \leq \frac{1}{2} h(\p)^2 \frac{ \max_{\theta \in [0, \alpha]}
		|\s^\T (\nabla\nabla r_i)(\xib+\theta \s)\s|}{\|\nabla \r(\xib) \s\|^2},
	\end{equation}
	which gives us the local error bound for a point-candidate pair:
	\begin{equation}
	|\Delta h(\xib, \s)| \le \frac{\sqrt{d_\Xi}}{2} h(\p)^2 \frac{ \displaystyle
		\max_{i=1,\ldots,d_\Xi}
		\max_{\theta \in [0, \alpha]}
		\left|\s^\T (\nabla\nabla r_i)(\xib+\theta \s)\s\right|}{\|\nabla \r(\xib) \s\|^2}.
	\end{equation}
	To estimate the error around $\p$ for all candidates,
	we can further bound~\eqref{eq:component-bound} more independently of $\s$ by using
	\begin{align}
	\|\nabla \r(\xib) \s\| &\ge \sigma_{d_\Xi}(\nabla \r(\xib)) > 0 \\
	\left|\s^\T (\nabla\nabla r_i)(\xib+\theta \s)\s\right| &
	\le \sigma_1((\nabla\nabla r_i)(\xib + \theta \s))
	\end{align}
	where $\sigma_{d_\Xi}$ is the smallest singular value of the Jacobian,
	which is positive as $\r$ is regular and
	$\sigma_1$ is the largest singular value of the Hessian.
	Following that, we can bound $R_i$ as
	\begin{align}
	|R_i(\xib, \vec s)| &\leq \frac{1}{2} h(\p)^2 \frac{ \max_{\b \zeta  \in [0, \alpha]}
		\sigma_1((\nabla\nabla r_i)(\xib + \theta \s))}{\sigma_{d_\Xi}(\nabla \r(\xib))^2} \\
	&\le  \frac{1}{2} h(\p)^2 \frac{ \displaystyle \max_{\b \zeta \in \bar B(\xib, \rho_{\xib})}
		\sigma_1((\nabla\nabla r_i)(\b \zeta))}{\sigma_{d_\Xi}(\nabla \r(\xib))^2}
	\label{eq:ball-ineq}
	\end{align}
	where
	\begin{equation}
	\rho_{\xib} = \frac{h(\p)}{\sigma_{d_\Xi}(\nabla \r(\xib))} \ge \alpha(\xib, \s)
	\end{equation}
	is the radius of the closed ball $\bar B(\xib, \rho_{\xib})$ centered at $\xib$. The
	inequality~\eqref{eq:ball-ineq} holds as the maximum is sought over a larger domain.
	This gives the local estimate around a point as
	\begin{equation} \label{eq:local-point}
	|\Delta h(\xib)| \le  \frac{\sqrt{d_\Xi}}{2} h(\p)^2 \frac{ \displaystyle
		\max_{i=1,\ldots,d_\Xi}
		\max_{\b \zeta \in \bar B(\xib, \rho_{\xib})}
		\sigma_1((\nabla\nabla r_i)(\b \zeta))}{\sigma_{d_\Xi}(\nabla \r(\xib))^2}.
	\end{equation}

	To obtain a simple global estimate, we take the maximum of~\eqref{eq:local-point} over $\xib$
	\begin{align}
	|\Delta h| &\le  \max_{\xib \in \Xi} \left( \frac{\sqrt{d_\Xi}}{2} h(\p)^2 \frac{
		\displaystyle \max_{i=1,\ldots,d_\Xi}
		\max_{\b \zeta \in \bar B(\xib, \rho_{\xib})}
		\sigma_1((\nabla\nabla r_i)(\b \zeta))}{\sigma_{d_\Xi}(\nabla \r(\xib))^2}\right) \\
	&\le  \frac{\sqrt{d_\Xi}}{2} h_M^2
	\frac{
		\displaystyle \max_{i=1,\ldots,d_\Xi} \max_{\xib \in \Xi}
		\max_{\b \zeta \in \bar B(\xib, \rho_{\xib})}
		\sigma_1((\nabla\nabla r_i)(\b \zeta))}{\min_{\xib \in \Xi} \sigma_{d_\Xi}^2(\nabla \r(\xib))}
	\\
	&\le \frac{\sqrt{d_\Xi}}{2} h_M^2 \frac{\sigma_{1,M} (\nabla\nabla\r)}{\sigma_{d_\Xi,
			m}^2(\nabla\r)},
	\end{align}
	where \begin{align}
	h_M &= \max_{\xib \in \Xi} h(\p),\\
	\label{eq:sing-max}
	\sigma_{1,M} (\nabla\nabla\r) &= \max_{i=1,\ldots,d_\Xi} \max_{\xib \in \Xi}
	\sigma_1((\nabla\nabla r_i)(\xib)), \\
	\sigma_{d_\Xi, m}(\nabla\r) &= \min_{\xib \in \Xi} \sigma_{d_\Xi}(\nabla \r(\xib)),
	\end{align}
	and the equality in~\eqref{eq:sing-max} holds since the inner maximum over local balls
	is superfluous when the maximum is sought over the whole domain.
	In particular, if $h$ is constant, the absolute error is of order $h^2$ and the
	relative error is linear.
\end{proof}

\bibliographystyle{siamplain}
\bibliography{references}
\end{document}